\documentclass[11pt,leqno]{amsart}

\usepackage[english]{babel}
\usepackage[all,2cell,emtex]{xy}
\usepackage{tikz-cd}
\tikzset{
  commutative diagrams/.cd,
  arrow style=tikz,
  diagrams={>=latex}}
\makeatletter
\tikzset{
  column sep/.code=\def\pgfmatrixcolumnsep{\pgf@matrix@xscale*(#1)},
  row sep/.code   =\def\pgfmatrixrowsep{\pgf@matrix@yscale*(#1)},
  matrix xscale/.code=%
    \pgfmathsetmacro\pgf@matrix@xscale{\pgf@matrix@xscale*(#1)},
  matrix yscale/.code=%
    \pgfmathsetmacro\pgf@matrix@yscale{\pgf@matrix@yscale*(#1)},
  matrix scale/.style={/tikz/matrix xscale={#1},/tikz/matrix yscale={#1}}}
\def\pgf@matrix@xscale{1}
\def\pgf@matrix@yscale{1}
\makeatother

\usepackage{enumitem}
\usepackage{url,amssymb,amsmath,bbm,xr,color,nicefrac,stmaryrd}
\usepackage[mathscr]{euscript}
\usepackage{mathrsfs}
\usepackage[OT2,OT1]{fontenc}
\usepackage[colorlinks=true,pagebackref=true]{hyperref} 

\input cyracc.def
\DeclareFontFamily{U}{russian}{}
\DeclareFontShape{U}{russian}{m}{n}
        { <5><6> wncyr5
        <7><8><9> wncyr7
        <10><10.95><12><14.4><17.28><20.74><24.88> wncyr10 }{}
\DeclareSymbolFont{Russian}{U}{russian}{m}{n}
\DeclareSymbolFontAlphabet{\mathcyr}{Russian}
\makeatletter
\let\@math@cyr\mathcyr
\renewcommand{\mathcyr}[1]{\@math@cyr{\cyracc #1}}
\makeatother

\newcommand{\modl}[1]{#1\textnormal{-}\mathrm{mod}}

\newcommand{\base}{\mathscr S}

\newcommand{\T}{\mathscr T} 
\newcommand{\C}{\mathscr C} 

\DeclareMathOperator{\uK}{\underline K}

\DeclareMathOperator{\uPic}{\underline{Pic}}
\DeclareMathOperator{\uZZ}{\underline{\ZZ}}

\DeclareMathOperator{\rk}{rk}

\DeclareMathOperator{\GW}{GW}

\newcommand{\E}{\mathbb E}
\newcommand{\F}{\mathbb F}
\newcommand{\un}{\mathbbm 1}

\newcommand{\MSL}{\mathbf{MSL}}
\newcommand{\KQ}{\mathbf{KQ}}
\newcommand{\KW}{\mathbf{KW}}

\newcommand{\KMW}{\mathbf K^\mathrm{MW}}

\newcommand{\uW}{\underline{\mathrm W}}
\newcommand{\HuW}{\mathbf H \underline{\mathrm W}}
\newcommand{\HM}{\mathbf H_{\mathrm{M}}}
\newcommand{\HMW}{\mathbf H_{\mathrm{MW}}}
\newcommand{\HW}{\mathbf H_{\mathrm{W}}}

\newcommand{\HH}{\mathrm H}
\newcommand{\HA}{\mathrm H_{\AA^1}}
\newcommand{\BA}{\mathrm H^{\AA^1}}

\newcommand{\GL}{\mathrm{GL}}
\newcommand{\SL}{\mathrm{SL}}

\newcommand{\thom}{\mathfrak t}

\newcommand{\vb}[1]{\langle#1\rangle}

\DeclareMathOperator{\codim}{codim}
\newcommand{\eff}{\mathrm{eff}}

\newcommand{\sm}{\mathrm{Sm}}

\DeclareMathOperator{\DM}{DM}
\DeclareMathOperator{\DMB}{DM_\mathcyr B}
\DeclareMathOperator{\SH}{SH}
\DeclareMathOperator{\DA}{D_{\AA^1}}
\DeclareMathOperator{\DAeff}{D_{\AA^1}^{\eff}}

\DeclareMathOperator{\DAminus}{D_{\AA^1,-}}

\DeclareMathOperator{\Id}{Id}

\DeclareMathOperator{\Hom}{Hom}
\DeclareMathOperator{\Maps}{Maps}

\DeclareMathOperator{\uHom}{\underline{Hom}} 
\DeclareMathOperator{\Pic}{Pic}
\DeclareMathOperator{\CH}{\mathrm{CH}}
\DeclareMathOperator{\wCH}{\widetilde{\mathrm{CH}}}
\DeclareMathOperator{\spec}{Spec}
\DeclareMathOperator{\Ind}{Ind}

\DeclareMathOperator{\ilim}{\varinjlim}
\DeclareMathOperator{\plim}{\varprojlim}

\newcommand{\twist}[1]{\langle #1 \rangle}
\DeclareMathOperator{\Th}{\mathrm{Th}}
\DeclareMathOperator{\Tw}{\mathrm{Tw}}
\DeclareMathOperator{\tdet}{\widetilde \det}
\DeclareMathOperator{\tomega}{\widetilde \omega}

\DeclareMathOperator{\K}{K}
\DeclareMathOperator{\Der}{D}
\DeclareMathOperator{\Spt}{Spt} 

\newcommand{\ZZ} {\mathbb Z}
\newcommand{\ZZot} {\mathbb Z[\nicefrac 1 2]}
\newcommand{\ot} {{\nicefrac 1 2}}
\newcommand{\QQ} {\mathbb Q}

\newcommand{\OO}{\mathcal O}

\renewcommand{\AA} {\mathbb A}
\newcommand{\PP} {\mathbb P}

\newcommand{\GG} {\mathbb{G}_m}

\newcommand{\GGx}[1] {\mathbb{G}_{m,#1}}

\newcommand{\nis}{\mathrm{Nis}}
\newcommand{\zar}{\mathrm{Zar}}

\title{On the rational motivic homotopy category}

\author{Fr\'ed\'eric D\'eglise}
\address{ENS de Lyon, UMPA, UMR 5669, 46 allée d'Italie 69364, Lyon Cedex 07, France}
\email{frederic.deglise@ens-lyon.fr}
\urladdr{http://perso.ens-lyon.fr/frederic.deglise/}
 
\author{Jean Fasel}
\address{Institut Fourier - UMR 5582, Universit\'e Grenoble-Alpes, CS 40700, 38058 Grenoble Cedex 9, France}
\email{jean.fasel@gmail.com}
\urladdr{https://www-fourier.univ-grenoble-alpes.fr/~faselj/}

\author{Adeel A. Khan}
\address{Institut des Hautes Études Scientifiques\\
35 route de Chartres\\
91440 Bures-sur-Yvette\\
France}
\email{khan@ihes.fr}
\urladdr{https://www.preschema.com/}

\author{Fangzhou Jin}
\address{School of Mathematical Sciences\\
Tongji University\\
Siping Road 1239\\
200092 Shanghai\\
China}
\email{fangzhoujin@tongji.edu.cn}
\urladdr{https://fangzhoujin.github.io/}

\date{\today}

\newtheorem{thm}{Theorem}[section]
\newtheorem{prop}[thm]{Proposition}
\newtheorem{lm}[thm]{Lemma}
\newtheorem{cor}[thm]{Corollary}
\newtheorem{thmi}{Theorem}
\newtheorem{cori}[thmi]{Corollary}
\newtheorem*{lmi}{Key Lemma}

\theoremstyle{remark} 
\newtheorem{rem}[thm]{Remark}
\newtheorem{notation}[thm]{Notation}
\newtheorem{ex}[thm]{Example}

\newtheorem{setup}[thm]{Setup}
\newtheorem{warn}[thm]{Warning}

\theoremstyle{definition} 
\newtheorem{df}[thm]{Definition}
\newtheorem{num}[thm]{}

\numberwithin{equation}{thm}

\begin{document}

\begin{abstract}
  We study the structure of the rational motivic stable homotopy category over general base schemes.
  Our first class of results concerns the six operations: we prove absolute purity, stability of constructible objects, and Grothendieck--Verdier duality for $\SH_\QQ$.
  Next, we prove that $\SH_\QQ$ is canonically $\SL$-oriented; we compare $\SH_\QQ$ with the category of rational Milnor--Witt motives; and we relate the rational bivariant $\AA^1$-theory to Chow--Witt groups.
  These results are derived from analogous statements for the minus part of $\SH[1/2]$.
\end{abstract}

\maketitle

\setcounter{tocdepth}{3}
\tableofcontents

\section{Introduction}

\subsection*{Summary of results}

Let $S$ be a scheme and $\SH(S)$ the stable motivic homotopy category.
This paper is dedicated to the study of its rationalization $\SH(S)_\QQ$.

The first relevant observation, due to Morel, is that after inverting $2$, $\SH(S)$ splits into two parts:
$$
  \SH(S)[1/2] \simeq \SH(S)_+ \times \SH(S)_-.
$$
In particular, one has a decomposition $\SH(S)_\QQ \simeq \SH(S)_{\QQ,+}\times \SH(S)_{\QQ,-}$.
The plus part is very well understood after the work of Cisinski--D\'eglise \cite{CD3}:

\begin{thmi}\label{thm:intro_SHQ+}
  The motivic $\infty$-category\footnote{%
    This is a refinement of the notion of motivic triangulated category from \cite{CD3}; see \cite{KhanThesis,KhanSix}.
  } $\SH_{\QQ,+}$ satisfies the following properties:

  \begin{enumerate}[label={\rm(\Roman*)}, leftmargin=*]
    \item\emph{Orientation.}\label{item:intro_SHQ+/orientation}
    It admits a canonical orientation.
    In particular, there are canonical isomorphisms
    $$
      \Th_X(v) \simeq \un_X\twist{n},
    $$
    for every virtual vector bundle $v$ of rank $n$ over a scheme $X$.
    Here we write $\twist{n} = (n)[2n]$ (Tate twist by $n$ and shift by $2n$) for an integer $n$.
    \item\emph{Absolute purity.}
    For any smoothable\footnotemark~morphism $f : X \to S$ of regular noetherian schemes, there is a canonical isomorphism
      $$ \un_{X}\twist{d} \simeq f^!(\un_{S}) $$
    in $\SH(X)_{\QQ,+}$, where $d = \rk(T_{X/S})$ is the relative virtual dimension of $f$.
    \footnotetext{A morphism is called \emph{smoothable} if it admits a global factorization through a closed immersion and a smooth morphism.}
    \item\emph{Finiteness.}
    Over quasi-excellent schemes, each of the six operations preserves constructible objects.
    \item\emph{Duality.}
    For every separated morphism of finite type $f : X \to S$ with $S$ quasi-excellent and regular, $f^!(\un_{S})$ is a dualizing object of $\SH(X)_{\QQ,+}$.
    \item\emph{Comparison.}\label{item:intro_SHQ+/comparison}
    There is a canonical equivalence of symmetric monoidal motivic $\infty$-categories
    $$
      \SH_{\QQ,+} \simeq \modl{\HM\QQ},
    $$
    where $\HM\QQ$ is the rational motivic cohomology spectrum.
    \item\emph{Bivariant theory.}
    Let $S$ be a regular noetherian scheme equipped with the dimension function $\delta=-\codim_S$.
    Then for every $S$-scheme $X$ of finite type, there are canonical isomorphisms
    $$
      \BA_{2n,n}(X/S, \QQ)_+
      := \big[ \un_X\twist{n}, f^!(\un_S) \big]_{\SH(X)_{\QQ,+}}
      \simeq \CH_{\delta=n}(X) \otimes_\ZZ \QQ
    $$
    for every integer $n$, where $f : X \to S$ is the structural map and $\CH_{\delta=n}(X)$ is the Chow group of algebraic cycles on $X$ of $\delta$-dimension $n$.
    \footnote{Note that the dimension convention here differs from the usual one for Chow groups (see for example \cite{Ful}), as the integer $\delta$ may be negative. For more details about our convention, see Definition~\ref{df:CHtilde} or \cite[Example 3.1.6]{BD1}.}
  \end{enumerate}
\end{thmi}

Our goal in this paper is to demonstrate analogous properties for the minus part $\SH_{\QQ,-}$ (and hence for $\SH_\QQ$ itself).
In fact, we are able to prove these just after $\ZZ[1/2]$-linearization (i.e., for the minus part of the $\AA^1$-derived category):

\begin{thmi}\label{thm:intro_SHminus}
  Consider the minus part of the $\AA^1$-derived category, i.e., the $\ZZ[1/2]$-linear motivic $\infty$-category
  $$
    \DAminus
    = \SH_- \otimes_{\Spt} {\Der(\ZZ)}
    \simeq (\SH \otimes_{\Spt} \Der(\ZZ))[1/2]_-,
  $$
  where $\Spt$ is the $\infty$-category of spectra, $\Der(\ZZ)$ the derived $\infty$-category of $\ZZ$-modules, and the tensor product is in the sense of \cite[\S 4.8]{LurieHA}.
  More generally, let $\DAminus(-, \Lambda) := \SH_- \otimes_{\Spt} \Der(\Lambda)$ for any commutative ring $\Lambda$.
  Then $\DAminus(-, \Lambda)$ satisfies the following properties:

  \begin{enumerate}[label={\rm(\Roman*)}, leftmargin=*]
    \item\emph{Weak orientation.}
    The motivic $\infty$-category $\DAminus(-, \Lambda)$ admits a canonical $\SL$-orientation.
    In particular, there are canonical isomorphisms
    $$
      \Th_X(v) \simeq \Tw_X(\tdet(v)),
    $$
    for every virtual vector bundle $v$ over a scheme $X$.
    Here $\tdet(v)$ denotes the graded determinant, and for a graded line bundle $(L,n)$, we write $\Tw_X(L,n) = \Th_X(L)\twist{n-1}$.
    \item\emph{Absolute purity.}
    For any smoothable morphism $f : X \to S$ of regular noetherian schemes, there is a canonical isomorphism
    $$
      \Tw_X(\tomega_{X/S})
      \simeq f^!(\un_{S})
    $$
    in $\SH(X)_-$, where $\tomega_{X/S}$ is the graded determinant of the virtual tangent bundle $T_{X/S}$ (see \ref{num:grdet}).
    \item\emph{Finiteness.}
    Over quasi-excellent schemes, each of the six operations preserves constructible objects.
    \item\emph{Duality.}
    Let $S$ be a quasi-excellent regular scheme and $K_S \in \SH(S)_-$ a $\otimes$-invertible object.
    Then for every separated morphism of finite type $f : X \to S$, $f^!(K_S)$ is a dualizing object of $\SH(X)_-$.
    \item\emph{Comparison.}\label{item:intro_SHminus/comparison}
    Assume $\Lambda$ is a subring of $\QQ$ with $2 \in \Lambda^\times$.
    There is a canonical equivalence of $\Lambda$-linear symmetric monoidal motivic $\infty$-categories
    $$
      \DAminus(-, \Lambda) \simeq \modl{\HW \Lambda},
    $$
    where $\HW \Lambda$ is the $\Lambda$-linear homotopical Witt spectrum.
    \item\emph{Bivariant theory.}
    Assume $\Lambda$ is a subring of $\QQ$ with $2 \in \Lambda^\times$.
    Let $S$ be a regular noetherian scheme equipped with the dimension function $\delta=-\codim_S$.
    For every $S$-scheme $X$ of finite type, there are canonical isomorphisms
    $$
      \BA_{2n,n}(X/S, \Lambda_-)
      := \big[ \un_X\twist{n}, f^!(\un_S) \big]_{\DAminus(X, \Lambda)}
      \simeq \wCH_{\delta=n}(X) \otimes_\ZZ \Lambda
    $$
    for every integer $n$, where $f : X \to S$ is the structural map and $\wCH_{\delta=n}(X)$ is the Chow--Witt group of quadratic cycles on $X$ of $\delta$-dimension $n$.
  \end{enumerate}
\end{thmi}

In fact, the motivic $\infty$-category $\SH_{-}$ also satisfies property \ref{item:intro_SHQ/purity}, when formulated correctly, as well as properties \ref{item:intro_SHQ/finiteness} and \ref{item:intro_SHQ/duality}.

Combining Theorems~\ref{thm:intro_SHQ+} and \ref{thm:intro_SHminus}, we get:

\begin{cori}\label{cor:intro_SHQ}
  The motivic $\infty$-category $\SH_{\QQ}$ satisfies the following properties (where (I)--(IV) are formulated exactly as in Theorem~\ref{thm:intro_SHminus}):

  \begin{enumerate}[label={\rm(\Roman*)}, leftmargin=*]
    \item\emph{Weak orientation.}\label{item:intro_SHQ/orient}
    \item\emph{Absolute purity.}\label{item:intro_SHQ/purity}
    \item\emph{Finiteness.}\label{item:intro_SHQ/finiteness}
    \item\emph{Duality.}\label{item:intro_SHQ/duality}
    \item\emph{Comparison.}\label{item:intro_SHQ/comparison}
    There is a canonical equivalence of $\QQ$-linear symmetric monoidal motivic $\infty$-categories
    $$
      \SH_{\QQ} \simeq \modl{\HMW\QQ},
    $$
    where $\HMW\QQ$ is the rational Milnor--Witt motivic cohomology spectrum.
    \item\emph{Bivariant theory.}\label{item:intro_SHQ/bivariant}
    Let $S$ be a regular noetherian scheme equipped with the dimension function $\delta=-\codim_S$.
    Then for every $S$-scheme $X$ of finite type, there are canonical isomorphisms
    $$
      \BA_{2n,n}(X/S, \QQ)
      := \big[ \un_X\twist{n}, f^!(\un_S) \big]_{\SH(X)_\QQ}
      \simeq \wCH_{\delta=n}(X) \otimes_\ZZ \QQ
    $$
    for every integer $n$, where $f : X \to S$ is the structural map and $\wCH_{\delta=n}(X)$ is the Chow--Witt group of quadratic cycles on $X$ of $\delta$-dimension $n$.
  \end{enumerate}
\end{cori}

We now discuss each of the properties (I)--(VI) in turn.

\subsection*{(I) Weak orientation}

Orientations are the main point of difference we encounter from the classical case of $\SH_{\QQ,+}$.
In fact, $\DAminus$ and $\SH_\QQ$ do not admit ($\GL$-)orientations in the sense of \cite[Def.~2.4.38]{CD3}.

Before explaining the weaker notion of $\SL$-orientation, we briefly recall the classical notion.
For any motivic $\infty$-category $\T$, there is a canonical map of presheaves of groupoids
$$
  \Th = \Th^\T : \uK \to \uPic(\T),
$$
sending a virtual vector bundle $v$ over a scheme $S$ to the associated \emph{Thom object}, a $\otimes$-invertible object $\Th_S(v) \in \T(S)$.
An \emph{orientation} of $\T$ amounts to the data of, for every virtual vector bundle $v$ of rank $n$ on a scheme $X$, a Thom isomorphism
$$
  \Th_X(v) \simeq \twist{n},
$$
where $\twist{n} = (n)[2n]$.%
\footnote{Recall that the \emph{Tate sphere} is defined as $\un_X(1):=\Th_X(\mathbb{A}^1_X)[-2]$; for any non-negative integer $n$, $\un_X(n):=\un_X(1)^{\otimes n}$, and if $n<0$, $\un_X(n):=\un_X(-n)^{\otimes (-1)}$.}
These isomorphisms are required to be functorial and compatible with exact triangles of perfect complexes.
Note that this is the same data as that of a factorization of $\Th$ through the rank map $\uK \to \ZZ$.

Now recall that an $\SL$-structure on a vector bundle is a trivialization of its determinant, and let $\uK^\SL$ denote the groupoid of $\SL$-oriented virtual vector bundles.
We define an \emph{$\SL$-orientation} of $\T$ as a factorization of the composite
$$
  \uK^\SL \to \uK \xrightarrow{\Th} \uPic(\T)
$$
through the rank map $\uK \to \ZZ$, where the first arrow is the forgetful map.
In other words, we have Thom isomorphisms only for vector bundles equipped with $\SL$-structures (and the isomorphism depends on the $\SL$-structure).
This notion is derived from the foundational work of Panin and Walter \cite{PaninWalter1}.

For \emph{any} virtual vector bundle $v$ of rank $n$, $v - \det(v)$ is a virtual vector bundle of rank $n-1$ with a canonical $\SL$-orientation.
Thus if $\T$ is $\SL$-oriented, we get canonical isomorphisms
$$
  \Th(v) \simeq \Tw(\det(v)),
$$
where we adopt the notation $\Tw(L) := \Th_X(L)\twist{n-1}$ for any graded line bundle $(L, n) \in \uPic$.

The fact that $\SH_{\QQ,+}$ is oriented (Theorem~\ref{thm:intro_SHQ+}\ref{item:intro_SHQ+/orientation}) is known by \cite[14.2.5 and Thm.~16.2.13]{CD3}.
Our result states that $\DAminus$ and $\SH_\QQ$ are canonically $\SL$-oriented (see Theorem~\ref{thm:DA_SL-oriented} and Corollary~\ref{cor:SHQ_SL-oriented}).
In particular, this means that every rational motivic spectrum $\E$ admits a canonical $\SL$-orientation.
This is a motivic analogue of the fact that any rational topological spectrum corresponds to an ``ordinary cohomology theory'', and in particular is complex-oriented (in the sense of Adams).

An important consequence concerns the \emph{fundamental class}
$$
  \eta_f : \Th_X(T_{X/S}) \to f^!(\un_S)
$$
associated to any smoothable lci morphism $f : X \to S$ in $\SH_-$ or $\SH_\QQ$ (see \cite[4.3.4]{DJK}).
Namely, this morphism simplifies to
$$
  \eta_f : \Tw_X(\tomega_{X/Y}) \to f^!(\un_S),
$$
where $\tomega_{X/S}$ is the graded determinant of the virtual tangent bundle $T_{X/S}$ (see~\ref{num:grdet} below). 
In particular, for a regular closed immersion $i : Z \to X$, the fundamental class determines a canonical element in twisted cohomology with support:
$$
  \eta_i
  \in H^{2c-2,c-1}_Z(X,\lambda_i)_\QQ
  = [X/X-Z,\Tw_X(\lambda_i,c)]_{\T(X)_\QQ},
$$
where $c$ is the codimension of $i$ and $\lambda_i$ the determinant of the normal bundle of $i$, which defines an element $(\lambda_i,c)\in \uPic$. 
In general, the above fundamental class lives in twisted $\AA^1$-bivariant theory, or rather its rational or minus part, and satisfies the classical properties of intersection theory, as explained in \cite{DJK} (compatibility with composition, transversal base change, excess intersection formula).
Moreover, it induces the usual wrong-way functorialities (covariance of twisted cohomology with respect to proper smoothable lci morphisms, contravariance of bivariant theories with respect to smoothable lci morphisms).

\subsection*{(II)--(IV) Absolute purity, finiteness and duality.}
These results are originally the main motivations of this paper.
Grothendieck's \emph{absolute (cohomological) purity conjecture} \cite[I 3.1.4]{SGA5} predicted that for a closed immersion $i:Y\to X$ between regular noetherian schemes of pure codimension $c$, and $n$ an integer invertible on $X$, there is an isomorphism $i^!(\mathbb{Z}/n)\simeq \mathbb{Z}/n(-c)[-2c]$ in the derived category of torsion \'etale sheaves on $Y$.
Though several cases were proved in the '60s (see \cite[I 5.1]{SGA5}, \cite[XVIII 3.2.5]{SGA4}), proving this property in the generality conjectured turned out to be a highly non-trivial problem: it was not until the '80s that a systematic approach for the general setting emerged, where a significant step was made by Thomason \cite{Thomason1}, who proved the conjecture assuming that all prime divisors of $n$ are sufficiently large, using the Atiyah-Hirzebruch spectral sequence for \'etale $K$-theory.
The general case was eventually solved by Gabber \cite{Fuji}, based on Thomason's method and the rigidity property for algebraic $K$-theory.
Further developments using a refinement of de Jong resolution of singularities, can be found in \cite[Exp.~XVI]{Gabber}.

Together with resolution of singularities, the absolute purity property lies at the very heart of Grothendieck's profound insight on the cohomological study of algebraic varieties. It has played a crucial role in the following applications:
\begin{enumerate}
\item The proof that the six functors on the derived category of torsion \'etale sheaves preserve constructible objects.

\item The proof of what is now called Grothendieck--Verdier local duality.  Most notably, if $S$ is a noetherian scheme equipped with a dimension function, there exists a \emph{dualizing complex} $K_S$ over $S$ such that $D_S=\uHom(-,K_S)$ is an auto-equivalence of categories. The functor $D_S$ then transforms $f_*$ (resp. $f^*$) into $f_!$ (resp. $f^!$).

\item The construction of Gysin morphisms and cycle classes, and the cohomological study of intersection theory.

\item The study of the coniveau spectral sequence.
\end{enumerate}
The attempt to answer these questions has inspired a number of important works based on new methods. Deligne managed to prove the first point unconditionally by avoiding the use of resolution of singularities \cite[Th. Finitude]{SGA4.5}. The second point was known in dimension at most one in \cite[Dualit\'e 1.4]{SGA4.5}, and proved by Gabber in general \cite[XVII]{Gabber}; furthermore, it implies the existence of the (self-dual) perverse $t$-structure over suitable base schemes, extending the fundamental work of \cite{BBD}.
The third point already appeared in \cite[XVIII]{SGA4}, and has since then been much worked on in the literature (see \cite{DJK}). The last point has also been discussed by a number of authors, such as \cite{BlochOgus} and \cite{CTHK}.

With the growing development of the theory of \emph{triangulated mixed motives} since the '90s, absolute purity naturally finds its place in this framework.
In fact, this property already appears implicitly in Beilinson's axiomatization modelled on the \'etale setting. Morel and Voevodsky's motivic homotopy theory came later \cite{MV}, together with a formalism of six functors in a similar fashion \cite{AyoubThesis}, in which absolute purity can be formulated for more general motivic categories, as highlighted in \cite[A.2]{CDEtale} and \cite[4.3]{DJK}.
Here are the first cases where absolute purity is known:
\begin{itemize}
\item For mixed motives with rational coefficients, it was first formulated and proved by Cisinski--D\'eglise in \cite[Thm.~14.4.1]{CD3}.
\item For modules over the algebraic $K$-theory spectrum, it is proven in \cite[Thm.~13.6.3]{CD3}.  Note that this serves as an ingredient in the proof of the previous case.
\item For \'etale motives, it is proven in \cite{AyoubEt} and \cite{CD5}.
\end{itemize}
These examples provide some evidence that such a property should hold in greater generality, and should philosophically be an addition to the six functors formalism: it was conjectured in \cite{Deg12} that absolute purity should also hold for the motivic homotopy category $\SH$ itself.

\bigskip

From the case of rational motivic cohomology, it follows from \cite[Thm.~16.2.13]{CD3} that absolute purity also holds for $\SH_{\QQ,+}$.
In this note, we show that it holds for $\SH_{-}$, and hence also for $\SH_{\QQ}$ (Theorem~\ref{thm:absolute_purity}).
These properties are known under the assumption (Resol) in \cite[2.4.1]{BD1}, but the difficulty is that the category $\SH_\QQ$ is not \emph{separated} in the sense of \cite[Definition 2.1.7]{CD3}.
The basis for our proof is the following simple but powerful observation (see Lemma~\ref{lm:fiber_Q_eq} and Example~\ref{ex:keylemma}(\ref{item:keylemma/SHminus})):

\begin{lmi}
  For every scheme $X$, there is a canonical equivalence of symmetric monoidal stable $\infty$-categories
  $$
    \SH(X)_- \simeq \SH(X_\QQ)_-,
  $$
  where $X_\QQ = X \times_{\spec{\ZZ}} \spec{\QQ}$, which commutes with the six operations.
\end{lmi}

This boils down to the calculation that $\SH(\spec{k})_- = 0$ for any field $k$ of positive characteristic.
The same holds for $\DAminus(-, \Lambda)$, for any commutative ring $\Lambda$, as well as the $\rho$-inverted variants (see Example~\ref{ex:keylemma}).

The Key~Lemma reduces the question of absolute purity to the equicharacteristic case, which is well-known (we include a proof in Appendix~\ref{apd:popescu}).
We thus deduce absolute purity for the rational, minus, and $\rho$-inverted parts of the motivic sphere spectrum.

The finiteness and duality results (properties \ref{item:intro_SHQ/finiteness} and \ref{item:intro_SHQ/duality}) are similarly deduced using the Key~Lemma from the equicharacteristic case, which is known by work of Bondarko--D\'eglise \cite{BD1} (see also \cite[Subsect.~3.1]{ElmantoKhan}).
See Propositions~\ref{prop:constructibility} and \ref{prop:duality}.

Finally, let us make the following remarks.
First, if $f : X \to S$ is a smoothable morphism of regular schemes, and $\tomega_{X/S}$ is the graded determinant of its relative virtual tangent bundle, then by \ref{item:intro_SHQ/purity} and \ref{item:intro_SHQ/duality}, the object
$$
  K_X = f^!(\un_S) \simeq \Tw_X(\tomega_{X/S})
$$
is a dualizing object in $\SH(X)_\QQ$.
Note the striking resemblance to Grothendieck--Serre duality for coherent cohomology \cite[Def.~III.2, Prop.~7.2]{HartshorneRD}.

\subsection*{(V) Comparison}

Let $\HM\QQ_S$ denote the rational motivic cohomology spectrum over a scheme $S$.
This can be defined either as the Beilinson motivic cohomology spectrum \cite[Def.~14.1.2]{CD3} or equivalently as the inverse image along the structural map $S \to \spec{\ZZ}$ of the rationalization of Spitzweck's integral motivic cohomology spectrum (by \cite[Thm.~7.14]{SpitzweckSpectrum}).
Modules over $\HM\QQ$ provide a good triangulated category of rational mixed motives, equivalent to Voevodsky's construction at least over excellent schemes \cite[Sect.~16.1]{CD3}.
The comparison equivalence (Theorem~\ref{thm:intro_SHQ+}\ref{item:intro_SHQ+/comparison})
$$
  \SH(S)_{\QQ,+} \simeq \modl{\HM\QQ}(S)
$$
was first observed by Morel over a field \cite{MorelSplitting} and generalized by Cisinski--D\'eglise \cite[Sect.~16.2]{CD3} to general bases.

For $\SH_{\QQ,-}$, an analogous statement was proven by Ananyevskiy--Levine--Panin over a field \cite[Sect.~4]{ALP}.
Our comparison generalizes this to arbitrary base schemes $S$ and to a $\Lambda$-linear statement for any commutative ring $\Lambda$.
Over any base scheme $S$, we define a motivic spectrum $\HW \Lambda_S$, the $\Lambda$-linear \emph{homotopical Witt spectrum} (see Definition~\ref{df:htp_Witt_singular}).
For $S$ regular, $\HW \Lambda_S$ is defined by the periodic $\GG$-spectrum
$$
  (\uW_S^\Lambda,\uW_S^\Lambda,\hdots)
$$
where $\uW_S$ is the Zariski sheaf associated to the Witt presheaf
$$
  S \mapsto W(S_\QQ) \otimes_\ZZ \Lambda
$$
on the site of smooth $S$-schemes of finite type (see Definition~\ref{df:unramified_Witt} and Paragraph~\ref{num:unram_Witt_spectrum}).

We then prove $\DA(S, \Lambda)_- \simeq \modl{\HW \Lambda}(S)$, hence in particular
$$
  \SH(S)_{\QQ,-} \simeq \modl{\HW \QQ}(S)
$$
over every scheme $S$ (see Corollary~\ref{cor:minus_DA&Witt_modl} and Paragraph~\ref{num:homotopical_Witt}).
Following \cite{ALP}, we can call the right-hand side the stable $\infty$-category of \emph{Witt motives} over $S$.
For $S$ regular, we deduce a computation of the minus part of the stable $\AA^1$-derived cohomology of $S$:
$$
  \HA^{n,i}(S,\Lambda)_-
  = \Hom_{DA(S,\Lambda)}(\un_S,\un_S(i)[n])_- \simeq H^{n-i}_\zar(S_\QQ,\uW) \otimes_\ZZ \Lambda.
$$

Finally we put the plus and minus parts together to get a model for $\SH_\QQ$.
Extending the original definition of \cite[Chap.~6]{BCDFO}, we define the rational \emph{Milnor--Witt motivic cohomology} spectrum over any scheme $S$ by
$$
  \HMW \QQ_S=\HM \QQ_S \oplus \HW \QQ_S,
$$
where $\HM \QQ_S$ is the rational motivic cohomology ring spectrum over $S$ (see Definition~\ref{df:HMWQ}).
We then deduce the comparison (Corollary~\ref{cor:rational_SH&DMW})
$$
  \SH(S)_{\QQ} \simeq \modl{\HMW^\QQ}(S).
$$
Thus we may view the right-hand side as the stable $\infty$-category of \emph{Milnor--Witt motives} over $S$.
For $S$ the spectrum of an infinite perfect field of characteristic unequal to $2$, this result was recently proven by Garkusha \cite{Garkusha}.

In particular, for a regular scheme $S$, we get the following computation
 of the rational stable homotopy classes of endomorphisms of the sphere spectrum (see Corollary~\ref{cor:rational_sphere&MW}):
$$
\big[\un_S,\un_S(i)[n]\big] \otimes_\ZZ \QQ \simeq K_{2i-n}^{(n)}(S)_\QQ \oplus H^{n-i}_\zar(S_\QQ,\uW)_\QQ.
$$
This result generalizes \cite[Remark 4.5]{ALP}. In the case where $S$ is the spectrum of a field (more generally when $S$ is a regular semi-local scheme over a field, see \cite[Theorem 10.12]{BachmannHoyois}), this recovers Morel's computation of the stable stems of $\mathbb{A}^1$-homotopy sheaves (\cite[Theorem 6.40]{MorelA1}).

We also apply these results to the dimensional homotopy t-structures of Bondarko--D\'eglise \cite{BD1} on $\SH_-$.
If $S$ is equipped with a dimension function $\delta$, then we consider the \emph{$\delta$-homotopy t-structure} on $\SH(S)_-$ and $\SH(S)_\QQ$ (Paragraph~\ref{num:recall_delta_htp}).
The Key~Lemma allows us to generalize the positivity/negativity criterion (Theorem~3 of \emph{op. cit.}) to this setting (Proposition~\ref{prop:BD331}).
Using this criterion we show that for a regular scheme $S$, the monoidal unit of $\DA(S,\Lambda)_-$ lies in the heart of the $\delta$-homotopy t-structure, when $\delta=-\codim_S$ (Corollary~\ref{cor:comput_HA1minus}).
(In contrast, this property does not hold for $\SH_{\QQ,+}$.)

\subsection*{(VI) Bivariant theory}

In any motivic $\infty$-category $\T$, the unit object represents a \emph{bivariant theory} by the formula
$$
  H^\T_{m,n}(X/S) := \big[ \un_X(n)[m], f^!(\un_S) \big]_{\T(X)}
$$
for $f : X \to S$ a separated morphism of finite type (see \cite[Def.~1.2.2]{DegliseBivariant}).
(Note that for $f=\Id_X$, we recover the cohomology theory represented by the unit, with opposite grading.)
The last part of our main result is concerned with identifying these groups explicitly when $\T$ is $\SH_{-}$ or $\SH_\QQ$, and $S$ is a regular base scheme equipped with the dimension function $\delta = -\codim_S$.

For $\SH_{\QQ,+}$ this was achieved by Bondarko--D\'eglise \cite[Theorem~2]{BD1}, using absolute purity and the niveau spectral sequence.
We follow the same ideas to get an analogous computation for $\SH_{-}$ and $\SH_\QQ$ (see Theorem~\ref{thm:biv_A1&ChowWitt}).
Let us make this computation explicit in the particular case $X=S$, a regular scheme.
Then we get an isomorphism of cohomology theories:
$$
  \HA^{2n,n}(X,\QQ)=\big[\un_X,\un_X\twist n\big]_\QQ \simeq \wCH^n(X) \otimes_\ZZ \QQ
$$
where the right-hand side is the group of ``rational equivalence classes of unramified quadratic cycles'' of codimension $n$ in $X$, i.e., the Chow--Witt group (see \cite{Fasel_wCH, FaselSrinivas}).
More explicitly, this is the middle cohomology group of a complex of the form:
$$
 \bigoplus_{\eta \in X^{(n-1)}} \KMW_1\big(\kappa_\eta,\omega_\eta\big)_\QQ
 \rightarrow \bigoplus_{x \in X^{(n)}} GW\big(\kappa_x,\omega_x\big)_\QQ
 \rightarrow \bigoplus_{s \in X^{(n+1)}} W\big(\kappa_s,\omega_s\big)_\QQ
$$
where for $x \in X$, $\kappa_x$ is the residue field 
 and $\omega_x$ the determinant of the cotangent complex of the immersion of $x$ in $X$,
 and $\K_1^{MW}$, $GW$, $W$ are respectively the twisted form of the
 first Milnor--Witt, Grothendieck--Witt, and Witt group of the given field
 (see paragraph \ref{num:RostSchmid} for more details).
 Note that when $2$ is invertible on $S$, this coincides with Theorem 4 of \cite{FaselSrinivas}.
 The above definition makes sense even without this assumption: the Key~Lemma tells us that the Witt part vanishes over the characteristic $2$ part of $S$!
 Note moreover that we give a formulation that also takes into account the presence of twists.

Using this computation we rule out the existence of Bondarko's Chow weight structure on $\SH(S)_\QQ$ for a large class of schemes $S$ (see Corollary~\ref{cor:chowwtnotexist}).
 Moreover, the computation is used in \cite{DF3} to obtain
 the so-called Borel character over arbitrary regular bases, which provides
 a rational isomorphism between hermitian K-theory and Milnor--Witt motivic cohomology extending the Chern character.

 \subsection*{Organization of the paper}

The paper is divided into 8 sections and contains 3 appendices.

The Key~Lemma is proved in Section~\ref{sec:key}.
 In Section~\ref{sec:abspur}, we prove properties \ref{item:intro_SHQ/purity}, \ref{item:intro_SHQ/finiteness}, and \ref{item:intro_SHQ/duality} of Theorem~\ref{thm:intro_SHminus} and Corollary~\ref{cor:intro_SHQ}.
 Section~\ref{sec:t-structure} deals with dimensional homotopy t-structures.
 We also obtain a useful description of the generators of $\SH_-$ and $\SH_\QQ$ (Proposition~\ref{prop:generation}).
 In Section~\ref{sec:W}, we give the construction of the \emph{homological Witt spectrum} and prove property \ref{item:intro_SHminus/comparison} for $\SH_-$.
 Section~\ref{sec:MW} deals with the same property \ref{item:intro_SHQ/comparison} for $\SH_\QQ$.
 In Section~\ref{sec:SL}, we reformulate Panin--Walter's theory of SL-orientations in terms of motivic categories, and explain its application to fundamental classes and Gysin morphisms.
 We also prove property~\ref{item:intro_SHQ/orient}, i.e., that $\SH_\QQ$ and $\SH_-$ are $\SL$-oriented.
 In Section~\ref{sec:bivariant}, we explain the computation of the $\delta$-niveau spectral sequence (a slight variation on the classical niveau spectral sequence using a dimension function $\delta$
 as introduced in \cite{BD1}) with respect to rational $\AA^1$-bivariant theory (Definition~\ref{df:bivariant_A1}), thus proving property~\ref{item:intro_SHQ/bivariant}.

In Appendix~\ref{apd:continuity}, we review the property of continuity in the language of motivic $\infty$-categories.
We use this property in Appendix~\ref{apd:efp} to extend the exceptional inverse image operation $f^!$ to essentially of finite presentation morphisms.
In Appendix~\ref{apd:popescu} we give a detailed proof of the absolute purity theorem for schemes defined over a field. (This result has already been noted in the literature, but without a detailed proof.)
 In Appendix~\ref{apd:hermitianK}, we prove absolute purity for hermitian K-theory,
 on the model of the case of Quillen K-theory \cite[Theorem 13.6.3]{CD3}.
 This result can be used to give another proof of Corollary~\ref{cor:intro_SHQ}\ref{item:intro_SHQ/purity} (at least for schemes over $\ZZ[1/2]$, see \cite{DFJKBorel}) and we include it for completeness.
 Finally in Appendix~\ref{apd:W}, we prove the analogue of the Key Lemma for the $2$-localised higher Witt groups of Balmer.
 We note that this readily implies the Gersten--Witt conjecture for arbitrary regular $\ZZot$-schemes,
 incidentally extending a result of Jacobson (see Remark~\ref{rem:Jacobson}).

\subsection*{Conventions}

All schemes are assumed quasi-compact and quasi-separated.
Regular schemes are assumed noetherian. 
Separated morphisms are assumed to be of finite type.
We will say that a morphism $f:X \rightarrow S$ of schemes is \emph{pro-\'etale} (resp. \emph{pro-smooth}) if $X$ is the projective limit of a cofiltered system of \'etale (resp. smooth) $S$-schemes with affine transition morphisms.
Given a property $\mathcal P$ of morphisms (e.g. \'etale, smooth, or separated), we will say that a morphism $f:X \rightarrow S$ of schemes \emph{essentially} satisfies $\mathcal P$ (or is essentially \'etale, smooth, or separated) if $X$ can be written as a projective limit of $S$-schemes satisfying $\mathcal P$ with affine and \'etale transition morphisms.
See Appendix~\ref{apd:efp} for details.

 We will freely use the language of \emph{motivic $\infty$-categories}, although the reader can mostly replace these with motivic triangulated categories \cite{CD3}.
 We will be more precisely concerned with the motivic stable homotopy category $\SH$ and the $\Lambda$-linear $\AA^1$-derived category $\DA(-,\Lambda)$ for a commutative ring of coefficients $\Lambda$.
 Recall that,  if $\T$ is any of these motivic $\infty$-categories,
 the motivic $\infty$-category $\T[1/2]$ obtained by inverting $2$ (or when $2 \in \Lambda^\times$) admits 
 following Morel a canonical orthogonal decomposition:
\begin{equation}\label{eq:+-decomposition}
\T[1/2] = \T_+ \times \T_-
\end{equation}
where $\T_+$ (resp. $\T_-$) is the (categorical) kernel of the projector
 $e_+=\frac{1-\epsilon} 2$ (resp. $e_-=\frac{1+\epsilon} 2$)
 --- here $\epsilon=-\langle -1\rangle$ in Morel's notation.
 We also note that $\T_-$ can be identified with $\T[1/2,\eta^{-1}]$, where $\eta$ is the algebraic Hopf map (see e.g. \cite[Rem.~3.6]{ALP}).

For a motivic $\infty$-category $\T$ and a scheme $S$, the full subcategory $\T_c(S)$ of \emph{constructible} objects is the thick stable subcategory generated by objects of the form $f_!f^!(\un_S)(n)$ for $f$ smooth separated of finite type and $n \in \ZZ$.
For the examples of $\T$ we consider, this will coincide with the subcategory of compact objects (see \cite[Prop.~1.4.11]{CD3}).

\bigskip

\noindent \textbf{Acknowledgments.}
F. D\'eglise received support from the French ``Investissements d'Avenir'' program, project ISITE-BFC (contract ANR-lS-IDEX-OOOB). F. Jin is partially supported by the DFG Priority Programme SPP 1786 Project ``Motivic filtrations over Dedekind domains'' and Marc Levine's ERC Advanced Grant QUADAG.
A. Khan was partially supported by SFB 1085 Higher Invariants, Universit\"at Regensburg, and from the Simons Collaboration on Homological Mirror Symmetry.

F. D\'eglise would like to thank Denis-Charles Cisinski whose influence via years of collaborations in motivic homotopy theory
 has been invaluable for the results obtained in this work.
F. Jin would like to thank Heng Xie for discussions about the paper \cite{Jacobson}.
A. Khan would like to thank Tom Bachmann, Marc Hoyois, and Charanya Ravi for helpful discussions related to this paper. We would like to thank the referee for carefully reading the paper and for the comments which helped improving the quality of the paper.

\section{The key lemma}
\label{sec:key}

\begin{num}\label{num:zero_part}
For a scheme $X$, denote by $X_\QQ = X\times_{\spec(\ZZ)}\spec(\QQ)$ the characteristic zero fiber.
We denote the inclusion by
\begin{equation*}
  \nu = \nu_X : X_\QQ \rightarrow X.
\end{equation*}
The following observation is the \textit{key lemma} of this paper:
\end{num}

\begin{lm}
\label{lm:fiber_Q_eq}
Let $\T$ be a continuous motivic $\infty$-category.
Suppose that for every field $k$ of positive characteristic, $\T(k)=0$.
Then for every scheme $X$, the functor
\begin{align}
  \nu_X^*:\T(X)\to\T(X_\QQ)
\end{align}
is an equivalence.
\end{lm}

\proof
If we write $X$ as a cofiltered limit of schemes $X_\alpha$ of finite type over $\spec{\ZZ}$ (with affine transition maps), then by continuity (Definition~\ref{df:continuity}) $\nu_X^*$ will be the colimit of the functors $\nu_{X_\alpha}^*$.
Thus we may assume $X$ is noetherian.
By Lemma~\ref{lm:pro-open immersion}, the co-unit $\nu^*\nu_* \to \Id$ is invertible since $\nu$ is a pro-open immersion.
It remains to show that the unit $\Id \to \nu_*\nu^*$ is also invertible.
By Proposition~\ref{prop:shriek conservative}, it will suffice to show that
\[
  i_x^! \to i_x^!\nu_*\nu^*
\]
is invertible for every point $x \in X$, where $i_x$ is the inclusion of the residue field at $x$ and where $i_x^!$ is defined as in Paragraph~\ref{num:f^! efp}.
By Propositions~\ref{prop:base change efp} and \ref{prop:Ex^*! efp} we may thus assume that $X$ is the spectrum of a field $k$.
If $k$ is of characteristic $0$, then $\nu$ is an isomorphism so the claim is obvious.
If $k$ is of positive characteristic, then both $\T(X)$ and $\T(X_\QQ)$ are zero by assumption, and the result follows.
\endproof

\begin{ex}\label{ex:keylemma}
  The condition of Lemma~\ref{lm:fiber_Q_eq} is satisfied for the following $\T$:
  \begin{enumerate}
    \item\label{item:keylemma/SHminus}
    $\T=\SH_-$, the minus part of the stable homotopy category.
    Indeed, recall that for any field $k$ there is a canonical isomorphism
      $$\Hom_{\SH(k)}(\un_k,\un_k) \simeq \GW(k),$$
    see \cite[Cor.~1.25]{MorelA1}\footnote{where the perfectness hypothesis is not necessary, see \cite[Thm.~10.12]{BachmannHoyois}}, where $\GW(k)$ is the Grothendieck--Witt ring of $k$.
    The decomposition $\SH[1/2] \simeq \SH_+ \times \SH_-$ corresponds to the decomposition
      $$ \GW(k)[1/2] \simeq \ZZ[1/2] \times W(k)[1/2], $$
    where $W(k)$ is the Witt ring (see \cite[Eqn.~(3.1)]{MorelA1}).
    Thus there is a canonical isomorphism
      $$ \Hom_{\SH(k)_-}(\un_{k,-},\un_{k,-}) \simeq W(k)[1/2]. $$
    But if $k$ is of positive characteristic and hence not formally real, then $W(k)[1/2]$ vanishes (see e.g. \cite[Subsect.~31.A]{EKM} or \cite[Chap.~VIII, Lem.~3.4]{LamQuadratic}).
    It follows that $\un_{k,-} \simeq 0$ and hence $\E_- = \E \otimes \un_{k,-} \simeq 0$ for every $\E \in \SH(k)$.

    \item\label{item:keylemma/DAminus}
    $\T=\DA(-, \Lambda)_-$, where $\Lambda$ is any commutative ring with $2 \in \Lambda^\times$.
    We may argue as in the first example, since the canonical map $\un_k \to \HA \ZZ_k$ induces an isomorphism on zeroth homotopy sheaves.
    Indeed, since $k$ is pro-smooth over its prime subfield, this follows by base change from the case where $k$ is perfect.
    In that case the claim is \cite[Th.~6.37]{MorelA1}.

    \item\label{item:keylemma/SHrhoinv}
    $\T=\SH(-)[1/\rho]$ or $\T=\DA(-)[1/\rho]$, where $\rho=-[-1]$ as in \cite[Corollary~19]{Bachmann}.
    This follows from the fact that, when $k$ is a field of positive characteristic, $\rho$ is nilpotent in $\SH(k)$ and $\DA(k)$, see the proof of \cite[Corollary~32]{Bachmann}.
  \end{enumerate}
  Moreover, each of these examples is continuous (see Example~\ref{ex:continuous}).
\end{ex}

\begin{rem}\label{rem:nu^* six operations}\leavevmode
\begin{enumerate}
\item
Under the assumptions of Lemma~\ref{lm:fiber_Q_eq}, we moreover have canonical identifications $\nu^*\simeq\nu^!$, where $\nu^!$ is defined as in Paragraph~\ref{num:f^! efp} (see Example~\ref{ex:f^! pro-etale}).
As this functor is an equivalence, its right adjoint $\nu_*$ is also a left adjoint.
Thus the functor $\nu^!$ admits a left adjoint $\nu_!$ (compare Warning~\ref{warn:f_! efp}).
Moreover, it immediately follows that both of the functors $\nu^* \simeq \nu^!$ and $\nu_* \simeq \nu_!$ commute with the six operations.
Similarly, symmetric monoidality of $\nu^*$ implies that $\nu_*\simeq\nu_!$ is also symmetric monoidal.\footnote{%
  By \cite[Prop.~2.5.5.1]{LurieSAG} (resp. its dual), $\nu_*$ is lax symmetric monoidal (resp. $\nu_!$ is colax symmetric monoidal).
  That is, $\nu_*\simeq\nu_!$ is both lax and colax symmetric monoidal, hence it is strict symmetric monoidal.
}

\item
Here is an alternative argument to show that the functor $\nu_X^*$ in Lemma~\ref{lm:fiber_Q_eq} is an equivalence for $\T=\SH(-)[1/\rho]$ (respectively $\T=\DA(-)[1/\rho]$): by \cite[Theorem 35]{Bachmann} there is a canonical equivalence $\SH(X)[1/\rho]\simeq\SH(X_r)$, where $X_r$ is the site associated to the \emph{real scheme} of $X$, which is a topological space whose points are pairs $(x,\sigma)$ where $x\in X$ is a point and $\sigma$ is an ordering of the residue field of $x$ (\cite[0.4.1]{Scheiderer}); since no field of positive characteristic is orderable, the morphism $\nu_X$ induces a homeomorphism (and therefore an equivalence of sites) $(X_\QQ)_r\simeq X_r$, which induces an equivalence $\SH(X)[1/\rho]\simeq\SH(X_\QQ)[1/\rho]$.
\end{enumerate}
\end{rem}

\section{Duality and absolute purity}
\label{sec:abspur}

In this section, we deduce from Lemma \ref{lm:fiber_Q_eq} several structural results for various motivic categories, corresponding to properties (II)--(IV) from the introduction.

\begin{num}
\label{num:Texample}
Let $\T$ be one of the following motivic $\infty$-categories:
\begin{itemize}
\item $\T=\SH_\QQ\simeq \DA(-,\QQ)$;
\item $\T=\SH_-$ or $\T=\DA(-,\Lambda)_-$, where $\Lambda$ is any commutative ring of coefficients;
\item $\T=\SH[1/\rho]$ or $\T=\DA[1/\rho]$, where $\rho$ is as in Example~\ref{ex:keylemma}\ref{item:keylemma/SHrhoinv}.
\end{itemize}
Note that in the last two cases $\T$ satisfies the condition in Lemma~\ref{lm:fiber_Q_eq}.
In the first case we may use the equivalence $\SH_\QQ\simeq\SH_{\QQ,+}\times\SH_{\QQ,-}$ to reduce to the plus and minus parts respectively. By \cite[Theorem~16.2.13]{CD3} the plus part is identified with Beilinson motives $\SH_{\QQ,+}\simeq\DM_{\textrm{\renewcommand\rmdefault{wncyr}\renewcommand\encodingdefault{OT2}\normalfont B},\QQ}$, for which the properties have already been established in \emph{op. cit}.
\end{num}

\begin{num}
The first application is the fact that the six functors preserve constructibility and what is called \emph{Grothendieck--Verdier local duality}, generalizing \cite[Theorem 2.4.9]{BD1}:
\end{num}
\begin{prop}[Constructibility]\label{prop:constructibility}
Let $\T$ be as in~\ref{num:Texample}.
 Given a morphism of schemes $f:X \rightarrow S$, the subcategory of constructible objects
 is stable under the following operations:
\begin{itemize}
\item $f^*$, $f_!$ for $f$ separated of finite type, $\otimes_S$;
\item $f_*$ for $f$ of finite type, $f^!$ for $f$ separated of finite type, $\uHom_S$, assuming in addition that $X$ and $S$ are quasi-excellent.
\end{itemize}
\end{prop}
\proof
If $\T$ satisfies the condition in Lemma~\ref{lm:fiber_Q_eq}, then by Remark~\ref{rem:nu^* six operations} one reduces to the case of characteristic zero schemes, for which the result follows from \cite[Theorem~2.4.9]{BD1}. For Beilinson motives the result is proved in \cite[Theorem~6.3.15]{CDEtale} (modulo the comparison of \cite[Theorem~16.1.2]{CD3}).
\endproof

\begin{prop}[Grothendieck--Verdier duality]\label{prop:duality}
Let $\T$ be as in~\ref{num:Texample}.
Let $S$ be a regular quasi-excellent scheme and $K_S \in \T(S)$ a $\otimes$-invertible object.
Then for any separated morphism $f:X\to S$ of finite type, the constructible object
  $$K_{X}=f^!(K_S)\in\T(X)$$
is \emph{dualizing}.
That is, for any constructible object $M \in \T(X)$, the canonical map
\begin{equation}\label{eq:can_loc_dual}
  M \rightarrow D_X(D_X(M))
\end{equation}
is invertible, where $D_X$ is the contravariant endofunctor of $\T(X)$ given by
  $$ D_X(M) = \uHom(M, K_X). $$
\end{prop}

\proof
If $\T$ satisfies the condition in Lemma~\ref{lm:fiber_Q_eq}, it suffices to show that the map~\eqref{eq:can_loc_dual} is an isomorphism after applying the equivalence $\nu^* : \T(X) \simeq \T(X_\QQ)$.
Since this equivalence commutes with the six operations (Remark~\ref{rem:nu^* six operations}), we are reduced to the equicharacteristic case, which is proven in \cite[Theorem~2.4.9]{BD1}. For Beilinson motives, the result is proven in \cite[Theorem~2.3.2]{CisLect}.
\endproof

\begin{num}
\label{num:abs_pur_def}
The most important consequence we obtain is the absolute purity property. Let $\T$ be a motivic $\infty$-category. Recall that an \emph{absolute object} $\E$ in $\T$ is a cartesian section, i.e. the data of an object $\E_X\in\T(X)$ for every scheme $X$ together with natural isomorphisms $f^*\E_Y\simeq \E_X$ for every morphism $f:X\to Y$ (compatible up to coherent homotopy). An absolute object $\E$ is said to satisfy \emph{absolute purity} if for any smoothable morphism $f : X \to Y$ between regular schemes with virtual tangent bundle $T_f$, the map
\begin{align}
\E_X\otimes\Th(T_f)\xrightarrow{}f^!\E_Y
\end{align}
induced by the fundamental class of $f$ (see \cite[Definition 4.3.7]{DJK}), is an isomorphism.
\end{num}
\begin{thm}[Absolute purity]\leavevmode
\label{thm:absolute_purity}
\begin{enumerate}
\item If $\T$ satisfies the condition in Lemma~\ref{lm:fiber_Q_eq}, then any absolute object in $\T$ satisfies absolute purity.
\item The minus part $\un_-$ and the rational part $\un_\QQ$ of the motivic sphere spectrum
 satisfy absolute purity in $\SH$.
\item Given any commutative ring $\Lambda$ with $2 \in \Lambda^\times$,
 the minus part of the $\AA^1$-derived $\Lambda$-linear ring spectrum $\Lambda_{-}$
 satisfies absolute purity in $\DA(-,\Lambda)$.
\end{enumerate}
\end{thm}
\begin{proof}
If $\T$ satisfies the condition in Lemma~\ref{lm:fiber_Q_eq}, then by Remark~\ref{rem:nu^* six operations} we are reduced to the characteristic zero case, where the result follows from Theorem~\ref{thm:abs_pur_eqchar}.
 Consider the second claim. For the ring spectrum $\un_-$, one applies the previous claim with $\T=\SH_-$.
 The case of the rational sphere spectrum $\un_\QQ \simeq \un_{\QQ,+}\oplus\un_{\QQ,-}$ follows from the case just treated and absolute purity for Beilinson motives (\cite[Theorem~14.4.1]{CD3}).
 The third claim follows from the first one with $\T=\DA(S,\Lambda)$.
\end{proof}

\section{The dimensional homotopy t-structure}
\label{sec:t-structure}

We next turn to applications about the dimensional homotopy t-structure \cite{BD1}.

\begin{num}\label{num:recall_delta_htp}
Recall the following definitions:
\begin{enumerate}[itemsep=0.2em]
\item
Let $\C$ be a compactly generated stable $\infty$-category which has small colimits. For $\mathcal{F}$ a family of objects in $\C$, there is a $t$-structure on $\C$, called the \emph{$t$-structure generated by $\mathcal{F}$}, where the subcategory of non-negative objects is the closure of $\mathcal{F}$ under extensions, positive shifts and small coproducts (\cite[Theorem~1.2.6]{BD1}, \cite[Proposition~2.1.70]{AyoubThesis}).

\item 
Let $\T$ be a motivic $\infty$-category. If $f:X\to S$ is a morphism of finite type, we denote $M^{BM}(X/S)=f_!\un_X\in\T(S)$ the \emph{Borel-Moore motive}. If $\E\in\T(S)$, we denote
  $$\E_{p,q}(X/S)=[M^{BM}(X/S)(q)[p],\E]_{\T(S)}$$
the \emph{bivariant group}.

\item
Let $\T$ be a motivic $\infty$-category and let $S$ be a scheme equipped with a dimension function $\delta$. The \emph{$\delta$-homotopy (perverse) $t$-structure} over $S$ is the $t$-structure $t_\delta$ on $\T(S)$ generated by the family of elements of the form $M^{BM}(X/S)(n)[n+\delta(X)]$ for $X$ a separated scheme of finite type over $S$ and $n\in\ZZ$ \cite[Definition 2.1.1]{BD1}.

\item
Let $S$ be a scheme and let $x:\spec{\kappa(x)}\to S$ be a topological point of $S$. An \emph{$S$-model} of $x$ is an affine regular $S$-scheme $X=\spec A$ of finite type, where $A$ is a subring of $\kappa(x)$ whose fraction field is $\kappa(x)$, such that $x$ coincides with the composition $\spec{\kappa(x)}\to X\to S$. A morphism of $S$-models is an $S$-morphism of the underlying schemes compatible with the structure of subring of $\kappa(x)$.
The category of $S$-models of $x$ is denoted $\mathcal{M}(x)$ \cite[Definition 3.2.1]{BD1}. If $S$ is excellent, $\mathcal{M}(x)$ is non-empty.

\item
Let $\E$ be an object in $\T(S)$. For $x$ a topological point of $S$ and $n,p\in \ZZ$, the \emph{fiber $\delta$-homology}
 of $\E$ in degree $p$ is defined as
\begin{align}
\hat H_p^\delta(\E)(x,n)=\varinjlim_{X\in\mathcal{M}(x)^{op}}\E_{2\delta(x)+p+n,\delta(x)+n}(X/S),
\end{align}
see \cite[Definition~3.2.3]{BD1}.\footnote{We remark that, according to the extension of the six functors
 formalism done in Appendix~\ref{apd:efp},
 we can write:
\begin{align*}
\hat H_p^\delta(\E)(x,n)=\E_{2\delta(x)+p+n,\delta(x)+n}(\kappa(x)/S)
 :=\Hom_{\T(x)}(\un_x(\delta(x)+n)[2\delta(x)+p+n],\nu_x^!\E).
\end{align*}
with $\nu_x:\spec(\kappa(x)) \rightarrow S$ the canonical immersion (Example~\ref{ex:residue efp}). In particular,  we could avoid
 the use of $S$-models from \cite{BD1}.
}
\end{enumerate}
Then we obtain the following generalization of \cite[Theorem 3.3.1]{BD1}:
\begin{prop}\label{prop:BD331}
Let $\T$ be as in~\ref{num:Texample} and let $S$ be an excellent noetherian scheme equipped with a dimension function $\delta$. Then for any object $\E\in\T(S)$ the following conditions are equivalent:
\begin{enumerate}
\item $\E$ is $t_\delta$-positive (resp. $\E$ is $t_\delta$-negative), that is, $\E\in\T(S)_{t_\delta\geqslant1}$ (resp. $\E\in\T(S)_{t_\delta\leqslant-1}$).

\item For any topological point $x$ of $S$ and integers $n\in\ZZ$ and $p\leqslant0$ (resp. $p\geqslant0$), $\hat H_p^\delta(\E)(x,n)=0$.
\end{enumerate}
\end{prop}
\proof
Suppose first that $\T$ satisfies the condition in Lemma~\ref{lm:fiber_Q_eq}, and let $\nu^*$ denote the equivalence of \emph{loc. cit}.
Denote by $\delta_\QQ$ the dimension function over $S_\QQ$ induced by $\delta$.
By \cite[Proposition~2.1.6]{BD1}, the functor $\nu^*$ is $t_\delta$-exact, so we find that $\E$ is $t_\delta$-positive if and only if $\nu^*\E$ is $t_{\delta_\QQ}$-positive. On the other hand, by \cite[Remark~3.2.5]{BD1} we know that $\hat H_p^\delta(\E)(x,n)=\hat H_p^{\delta_\QQ}(\nu^*\E)(x,n)$ if $x$ is a point of $S_\QQ$, and if $x$ is not a point of $S_\QQ$ we have $\hat H_p^\delta(\E)(x,n)=0$ by the condition in Lemma~\ref{lm:fiber_Q_eq}. Therefore we are reduced to prove the same statement for $\nu^*\E\in\T(S_\QQ)$, which follows from \cite[Example~3.2.14 and Theorem~3.3.1]{BD1}.

For Beilinson motives the result already follows from \cite[Theorem 3.3.1]{BD1}, since $\DM_{\textrm{\renewcommand\rmdefault{wncyr}\renewcommand\encodingdefault{OT2}\normalfont B},\QQ}$ is separated and therefore satisfies the condition (Resol) in \cite[2.4.1]{BD1}.
Hence we get the claim also for $\T = \SH_\QQ$.
\endproof
We deduce the analogues of the properties that are consequences of \cite[Theorem 3.3.1]{BD1}, as stated in \cite[Corollaries 3.3.5, 3.3.6, 3.3.7, 3.3.9, 3.3.10]{BD1}. In particular, the $t$-structure $t_\delta$ is non-degenerate.
\end{num}

\begin{num}
Over the ring of integers in a number field, we also obtain the following description of generators analogous to \cite[Corollary 2.4.8]{BD1}:
 
\begin{prop}[Generation]
\label{prop:generation}
Let $\T$ be as in~\ref{num:Texample}, and let $S$ be one of the following schemes:
\begin{enumerate}
\item \label{num:rsfirstcase}
$S$ is the spectrum of a commutative ring which is either
\begin{itemize}
\item the localization of the ring of integers of a number field at a set of primes,
\item or the ring of integers in a local field.
\end{itemize}
\item
$S$ is a scheme such that $-1$ is a sum of squares in every residue field.
\end{enumerate}
Then the stable $\infty$-category $\T(S)$ coincides with the localizing subcategory generated by objects of the form $f_*\un_X(n)$, where $f:X\to S$ is a projective morphism with $X$ a regular scheme and $n\in\mathbb{Z}$.
\end{prop}

\proof
For Beilinson motives the result is proved in \cite[Theorem 15.2.3]{CD3}.
Thus if $-1$ is a sum of squares in every residue field of $S$, then by \cite[Corollary~16.2.14]{CD3} we have $\SH(S,\QQ)\simeq\DM_{\textrm{\renewcommand\rmdefault{wncyr}\renewcommand\encodingdefault{OT2}\normalfont B}}(S,\QQ)$, and the result follows.

Thus we assume $S=\spec(R)$ is as in~\eqref{num:rsfirstcase}.
If $R$ has equal characteristic, the result is already proved in \cite[Corollary~2.4.8]{BD1}.
Otherwise, we have $S_\QQ=\spec{K}$ where $K$ is the fraction field of $R$.
If $\T$ satisfies the condition in Lemma~\ref{lm:fiber_Q_eq}, we only need to show that $\T(\spec{K})$ is generated as a localizing subcategory by objects of the form $\nu^* f_*\un_X(n) = \nu^* M^{BM}(X/S)(n)$, where $f:X\to S$ is a projective morphism with $X$ regular and $n\in\mathbb{Z}$.
By \cite[Corollary 2.4.8]{BD1} and noetherian induction, it suffices to show that for any irreducible smooth projective $K$-scheme $X$, there exists a non-empty subscheme $U \subseteq X$ such that $M^{BM}(U/K)\in\nu^*(\T')$. By de Jong's theory of alterations (see \cite[1.2.2]{Scholl}), there exists a projective generically finite surjection $X'\to X$ such that $X'$ has a projective regular model over $S$. By the proof of \cite[Corollary~2.4.8]{BD1}, there exists a non-empty open subscheme $U \subseteq X$ such that $M^{BM}(U/K)$ is a direct summand of an element in $\nu^*(\T')$, and the result follows.
\endproof

Applying Proposition~\ref{prop:generation}, we have the following analogue of \cite[Theorem 2.4.3]{BD1}:
\begin{prop}
\label{prop:BD243}
Let $\T$ be as in~\ref{num:Texample} and let $S$ be as in Proposition~\ref{prop:generation}. Then the homotopy $t$-structure on $\T$ agrees with the $t$-structure generated by elements of the form $M^{BM}(X/S)(n)[n+\delta(X)]$ for $X$ a regular scheme which is projective over $S$ and $n\in\ZZ$.
\end{prop}
Note that via Lemma~\ref{lm:fiber_Q_eq}, we have avoided using Proposition~\ref{prop:BD243} in the proof of Proposition~\ref{prop:BD331}, therefore removing the restriction on $S$.
\end{num}

\begin{num}
\label{num:eff_cat}
From Remark~\ref{rem:nu^* six operations} it follows that the functor $\nu^*$ preserves the $\delta$-effective subcategory (\cite[Definition~2.2.1]{BD1}). Consequently the analogues of Propositions~\ref{prop:BD331} and \ref{prop:BD243} for the effective category also hold, by a similar argument.
\end{num}

\section{The minus \texorpdfstring{$\AA^1$}{A1}-derived category and Witt motives}
\label{sec:W}

\begin{num}\label{num:symm_bdl}
We recall the classical definition of the Witt ring $W(S)$ of an arbitrary scheme $S$,
 following Knebusch \cite[I, \S 5]{Kne}.
 Given a scheme $S$, a \emph{symmetric bundle} (called ``space'' in \emph{loc. cit.})
 over $S$ is the data $(V,\phi)$ of a vector bundle $V$ over $S$
 together with a non-degenerate symmetric bilinear form $\phi:V \otimes V \rightarrow \AA^1_S$;
 \emph{i.e.} $\phi$ is a morphism of vector bundles such that $\phi \circ \epsilon=\phi$, where $\epsilon$ is the
 isomorphism exchanging factors, and the adjoint map $\phi':V \rightarrow \uHom(V,\AA^1_S)$ is an isomorphism.
 The notion of morphism of symmetric bundles is clear, and an isomorphism is called an \emph{isometry}.

Then $W(S)$ is the monoid of isometry classes of symmetric bundles, modulo the sub-monoid consisting of classes of metabolic symmetric bundles (see p.~133 and p.~145 in \cite{Kne}). Tensor product of symmetric bundles induces a product on the abelian group $W(S)$, which turns it into a commutative ring.
 It is clear from the definition that $W(S)$ is contravariantly functorial in $S$, with respect to arbitrary morphisms of schemes, so that $W$ defines a presheaf of commutative rings on the category of all schemes.
\end{num}
\begin{df}\label{df:unramified_Witt}
For any scheme $S$, we define the \emph{Witt sheaf} $\uW_S$ on the site of smooth $S$-schemes as the Zariski sheaf associated with the presheaf $X \mapsto W(X)$.
Given a subring $\Lambda$ of $\QQ$, we set $\uW_S^\Lambda=\uW_S \otimes_\ZZ \Lambda$.
\end{df}

\begin{rem}\label{rem:terminology_unramified_Witt}
According to the Gersten conjecture for the Witt ring
(known in the equicharacteristic case by \cite{BalmerGillePaninWalter}),
 it is well-known that for a regular scheme $S$ over a field,
 the Witt sheaf $\uW_S$ agrees with the unramified Witt sheaf, i.e., the sheaf of unramified elements in the Witt ring of the fraction field (in other words, the $0$-th cohomology of the Gersten--Witt complex).
 In fact, according to Corollary~\ref{cor:Gersten--Witt_conj},
 we further obtain that for \emph{any} regular scheme $S$,
 the $\Lambda$-linear Witt sheaf $\uW_S^\Lambda$ agrees with the $\Lambda$-linear unramified Witt sheaf, as long as $2\in \Lambda^\times$.
\end{rem}

\begin{rem}
Note that $\uW_S$ can also be described as the restriction to $\sm_S$ of the sheaf $\uW$ associated with $W$ on the big Zariski site of schemes of finite presentation over $S$.
 In particular, for any morphism of schemes $f:X \rightarrow S$ one gets a morphism of sheaves:
\begin{equation}\label{eq:funct1_uW}
\uW_S \rightarrow f_*(\uW_X),
\end{equation}
or equivalently by adjunction:
\begin{equation}\label{eq:funct2_uW}
f^*(\uW_S) \rightarrow \uW_X.
\end{equation}
\end{rem}

\begin{rem}\label{rem:W_S singular}
The above definition will be meaningful for us only when $S$ is regular.
As will be clear from the next results, this is explained by the failure of the Gersten--Witt conjecture over singular schemes.
\end{rem}

We will use the following basic lemma on Witt groups,
 for which we could not find a reference in the literature.
\begin{lm}\label{lm:Witt_lim}
Let $X$ be the limit of a cofiltered system of schemes $(X_i)_{i \in I}$ with affine transition maps. Then the canonical maps:
\begin{align*}
\ilim_{i \in I} W(X_i) & \rightarrow W(X)
\end{align*}
are isomorphisms.
Moreover, for any integer $n \geqslant 0$, the canonical maps
\begin{equation}\label{eq:unramf_Witt_coh&lim}
\ilim_{i \in I} H^n(X_i,\uW_{X_i}) \rightarrow H^n(X,\uW_X)
\end{equation}
are isomorphisms.
\end{lm}
\begin{proof}
Recall that, under these assumptions, the category of vector bundles over $X$ is equivalent to the 2-colimit of the category of vector bundles over $X_i$ (\cite[8.5.3--8.5.5]{EGA4}).
This implies the same result for symmetric bundles, by definition of the latter (Paragraph~\ref{num:symm_bdl}).
Moreover, if $V$ is a symmetric bundle over $X$, inverse image of a symmetric bundle $V_i$ over $X_i$ for some $i\in I$, then it is metabolic if and only if $V_j$ is metabolic for some $j\ge i$.
Thus the claim follows by definition of the Witt group.
The second claim then follows using the fact that the Zariski topology is coherent.
\end{proof}

\begin{thm}\label{thm:representability_Witt}
Let $S$ be a regular scheme over a prime field $k$ of characteristic not $2$.
 Then the following assertions hold:
\begin{enumerate}
\item
\label{num:WittZarNis}
 The Zariski sheaf $\uW_S$ is a Nisnevich sheaf and for any smooth $S$-scheme $X$,
 the following comparison map is an isomorphism:
\begin{align}
H^*_\zar(X,\uW_S) \simeq H^*_\nis(X,\uW_S).
\end{align}
Moreover, $\uW_S$ is strictly $\AA^1$-local\footnote{meaning that 
 its Nisnevich cohomology presheaves over $\sm_S$ are $\AA^1$-homotopy invariant.}.
\item 
\label{num:Hopfiso}
The algebraic Hopf map $\eta:(\GGx S,1) \rightarrow S_+$
 induces an isomorphism:
\begin{align}
\label{eq:Witt_Susp}
\eta_*:(\uW_S)_{-1}=\uHom(\Lambda(\GG,1),\uW_S) \simeq \uW_S,
\end{align}
where the internal Hom is taken in the effective $\AA^1$-derived category.
\end{enumerate}
\end{thm}
\begin{proof}
Consider assertion~\eqref{num:WittZarNis}.
When $S$ is smooth over $k$, the claim follows from \cite{Panin}.
 The case of an arbitrary regular scheme $S$
 can be obtained using Popescu's theorem (Theorem~\ref{thm:popescu}), the preceding lemma, and the fact that the smooth Zariski and Nisnevich toposes over $S$ are coherent (\cite[Exp.~VI, 2.3]{SGA4}).
 Alternatively, the reader can check that the arguments of \cite{Panin} actually work under our more general assumptions.

Then, the assertion~\eqref{num:Hopfiso} is a consequence of the localization exact sequence for triangular Witt groups,
 given the isomorphism $\Lambda(\GGx S,1) \simeq \Lambda(\Th(\AA^1_S))[-1]$
 in $\DAeff(S,\Lambda)$.\footnote{We refer the interested
 reader to the arguments of \cite[Lemma~3.8]{CF1} for details.}
\end{proof}

\begin{num}\label{num:unram_Witt_spectrum}
Let us now consider a regular $\QQ$-scheme $S$ and $\Lambda \subseteq \QQ$ a subring such that $2 \in \Lambda^\times$.
Recall that $\DAeff(S,\Lambda)$ can be seen as the full subcategory
 of the derived $\infty$-category of $\Lambda$-linear Nisnevich sheaves over $\sm_S$ spanned by $\AA^1$-local complexes.
 In particular, under the assumption of the previous theorem,
 the Witt sheaf $\uW^\Lambda_S$ defines such a complex
 (concentrated in degree $0$). Moreover,
 according to Theorem \ref{thm:representability_Witt}~\eqref{num:Hopfiso}, it is an infinite loop space with respect to $\un_S(1)[1]$.
 More precisely, there is an object $\HuW^\Lambda_S \in \DA(S,\Lambda)$ defined by the $\GG$-spectrum given by the collection
$$
(\uW^\Lambda_S,\uW^\Lambda_S,\hdots)
$$
together with the isomorphisms~\eqref{eq:Witt_Susp}.
 We call $\HuW^\Lambda_S$ the \emph{$\Lambda$-linear homotopical Witt spectrum} over $S$ (see Corollary~\ref{cor:comput_HA1minus} for the justification of the adjective ``homotopical'').
 By construction, this spectrum is $\eta$-periodic and satisfies:
$$
\Hom_{\DA(S,\Lambda)}(\Lambda_S(X),\HuW^\Lambda_S(i)[n])=H^{n-i}_\zar(X,\uW_S) \otimes_\ZZ \Lambda.
$$
The right-hand side is also the $(n-i)$-th cohomology group of the $\Lambda$-linearized Gersten--Witt complex.
 The commutative monoid structure on $\uW^\Lambda_S$ extends
 to a commutative ring structure on the spectrum $\HuW^\Lambda_S$,
 in such a way that the above isomorphism is compatible with products.

Note that it follows from Theorem \ref{thm:representability_Witt}~\eqref{num:Hopfiso}
that $\HuW^\Lambda_{S,+}=0$ or equivalently, that $\HuW^\Lambda_{S,-}\simeq\HuW^\Lambda_S$
 is an object of $\DA(S,\Lambda)_-$.
\end{num}

The Key Lemma~\ref{lm:fiber_Q_eq} allows us to define $\HuW^\Lambda_S$ in the mixed characteristic case:
\begin{df}\label{df:HW}
Let $S$ be an arbitrary regular scheme and $\Lambda \subseteq \QQ$
 a subring such that $2 \in \Lambda^\times$.
We define the $\Lambda$-linear homotopical Witt ring spectrum
 over $S$ as
$$
  \HuW^\Lambda_S := \nu_*\big(\HuW^\Lambda_{S_\QQ}\big),
$$
where $\nu:S_\QQ \rightarrow S$ is the inclusion of the characteristic zero fiber as in Paragraph~\ref{num:zero_part}.
\end{df}

\begin{rem}\label{rem:justif_df_uW}
While this definition may seem a bit artificial, the main point is that it allows us to bypass a missing part of the theory of higher Witt groups,
 which at the moment is not well-behaved over $\mathbb F_2$-schemes
 (for example, the localization exact sequence is not known for such schemes).
 Let us first mention that a recent work, \cite{realK}, may correct the situation.
 Moreover, Proposition~\ref{prop:Witt_2inverted} shows that
 one can extend the definition of $\HuW_S^\Lambda$ in Paragraph~\ref{num:unram_Witt_spectrum} 
 to the case of $S$ a regular $\ZZot$-scheme (with $\Lambda$ still a subring of $\QQ$ such that $2 \in \Lambda^\times$).
 In particular, in this latter case, the ring spectrum $\HuW_S^\Lambda$
 represents the cohomology of the $\Lambda$-linear Gersten--Witt complex
 (see Corollary~\ref{cor:Gersten--Witt_conj}).
\end{rem}

\begin{num}
Before moving to our main theorem, let us prove an intermediate lemma
 (which will be strengthened in Remark~\ref{rem:htp_Witt_singular}).
 Let $f:X \rightarrow S$ be an arbitrary morphism of regular schemes.
 Then there is a canonical morphism of ring spectra
$$
 \Theta_f:f^*\big(\HuW^\Lambda_S\big) \to \HuW^\Lambda_X
$$
in $\DA(X,\Lambda)$, induced by applying $\nu_*$ to the canonical map of $\GG$-spectra over $X_\QQ$
$$
\big(f^*\uW^\Lambda_{S_\QQ},f^*\uW^\Lambda_{S_\QQ},\hdots\big)
 \rightarrow \big(\uW^\Lambda_{X_\QQ},\uW^\Lambda_{X_\QQ},\hdots\big),
$$
which in turn comes from the functoriality~\eqref{eq:funct2_uW}.
\end{num}
\begin{lm}\label{lm:W commutes with p^*}
With notation as above, assume moreover that $X$ is pro-smooth over $S$.
 Then $\Theta_f$ is an isomorphism.
\end{lm}
\begin{proof}
By definition, we may assume that $S$ is defined over $\QQ$.
Since $f$ is pro-smooth, the functor $f^*$ preserves $\AA^1$-locality of complexes of $\Lambda$-linear Nisnevich sheaves (this follows using continuity, see Appendix~\ref{apd:continuity}, from the smooth case).
Thus it suffices to note that the map~\eqref{eq:funct2_uW} is an isomorphism
 of Nisnevich sheaves over $\sm_X$.
\end{proof}

\begin{thm}\label{thm:Witt motives}
Let $S$ be a regular scheme and $\Lambda \subseteq \QQ$ with $2 \in \Lambda^\times$.
 Then the unit map of the ring spectrum $\HuW_S$ induces an isomorphism of ring spectra
\begin{equation}\label{eq:W unit}
\Lambda_{S,-} \rightarrow \HuW^\Lambda_S
\end{equation}
in $\DA(S, \Lambda)$.
\end{thm}
\begin{proof}
Again using the Key~Lemma~\ref{lm:fiber_Q_eq}, and the fact that $\nu_*(\Lambda_{S_\QQ,-}) \simeq \Lambda_{S,-}$ (by monoidality, see Remark~\ref{rem:nu^* six operations}), we are reduced to the case where $S$ is a regular $\QQ$-scheme.
By Zariski descent \cite[2.3.8]{CD3} and Lemma~\ref{lm:W commutes with p^*}, we may also assume that $S$ is affine.
By Popescu's desingularization theorem (Theorem~\ref{thm:popescu}), any such $S$ is pro-smooth over $\QQ$, so by Lemma~\ref{lm:W commutes with p^*} we reduce to $S=\spec{\QQ}$.
Then the result follows from Proposition~41 of \cite{Bachmann}.
\end{proof}

\begin{cor}\label{cor:comput_HA1minus}
Let $S$ be a regular scheme and $\Lambda \subseteq \QQ$ with $2 \in \Lambda^\times$.
Then we have:
\begin{enumerate}
\item
There are canonical isomorphisms
$$
H^{n,m}_{\AA^1}(S,\Lambda)_- \simeq H^{n-m}_\zar(S_\QQ,\uW) \otimes_\ZZ \Lambda
$$
for any integers $(n,m) \in \ZZ^2$, compatible with pullbacks and products.

\item
The ring spectrum $\Lambda_{S,-} \simeq \HuW_S^\Lambda$ is in the heart of the $\delta$-homotopy t-structure on $\DA(S,\Lambda)$ for $\delta=-\codim_S$ (see Paragraph~\ref{num:recall_delta_htp}).

\item
The canonical map $\KW_S^\Lambda \rightarrow \HuW^\Lambda_S$ induces an isomorphism
$$
H_0^\delta\big(\KW_S^\Lambda\big) \rightarrow \HuW^\Lambda_S
$$
where $\KW_S^\Lambda$ is the spectrum representing the $\Lambda$-linearized higher Witt groups.
\end{enumerate}
\end{cor}

\begin{proof}
The first isomorphism is immediate from Theorem~\ref{thm:Witt motives}. To get the second assertion, one uses the following computation
 given a point $x \in S$:
$$
\hat H_p^\delta(\Lambda_{S,-})(x,n) \simeq H^{n-p,n}_{\AA^1}(\kappa(x),\Lambda)_-
 \simeq H^{-p}_\zar(\kappa(x),\uW) \otimes_\ZZ \Lambda
$$
which is obtained using absolute purity for $\un_-$ (Theorem~\ref{thm:absolute_purity}),
 as in \cite[Ex. 3.2.8]{BD1}.\footnote{Note the isomorphism depends on a choice
 of a trivialization of the normal bundle of $x$ in $X_{(x)}$.}
 Then Proposition~\ref{prop:BD331} gives the assertion.
 Finally, the last assertion can be reduced to the case of (characteristic zero)
 fields. Then the $\delta$-homotopy t-structure is just Morel's homotopy t-structure
 (\cite[Ex. 2.3.5(2)]{BD1}), and the computation follows from \cite[Th. 1, Rem. 1]{ALP}.
 \end{proof}

\begin{rem}
Note further that one obtains:
  $$\hat H_0^\delta(\Lambda_{S,-})(x,n) \simeq W\big(\kappa(x)\big)_\Lambda.$$
This isomorphism
 can be made canonical after tensoring with $\det\big(N_x(X_{(x)})\big)^\times$.
\end{rem}

\begin{cor}\label{cor:minus_DA&Witt_modl}
Let $S$ be a regular scheme and $\Lambda \subseteq \QQ$ with $2 \in \Lambda^\times$.
The canonical functor
$$
\DA(S,\Lambda)_- \rightarrow \DA(S,\Lambda)
$$
induces an equivalence of symmetric monoidal $\infty$-categories
$$
\DA(S,\Lambda)_- \simeq \modl{\HuW_S^\Lambda}
$$
where the right hand-side is the $\infty$-category of modules
 over $\HuW_S^\Lambda$, seen as a commutative algebra in the
 monoidal $\infty$-category $\DA(S,\Lambda)$.

Moreover, the $\delta$-homotopy $t$-structure
 on $\DA(S,\Lambda)_-$
 (Paragraph~\ref{num:recall_delta_htp})
 is non-degenerate and the unit lies in the heart
 (for $\delta=-codim_S$).
\end{cor}

\begin{rem}
The objects of the category $\modl{\HuW_S^\Lambda}$ are called \emph{Witt motives}
 in \cite[Section 4]{ALP}. One can obtain a somewhat more concrete description
 of this category along the lines of \cite{ALP}.
 Consider the symmetric monoidal Grothendieck abelian category
 $\mathrm{Sh}_{\uW}(S)$ of modules over the sheaf $\uW_S^\Lambda$ in the category of Nisnevich sheaves.
 Then define the $\Lambda$-linear $\infty$-category of \emph{Witt motives}
 $\DM_W(S,\Lambda)$ as the $\AA^1$-localization and $\PP^1$-stabilization of the derived $\infty$-category $\Der\big(\mathrm{Sh}_{\uW}(S)\big)$.
 The proof of the identification $\DM_W(S,\Lambda) \simeq \modl{\HuW_S^\Lambda}$ is analogous to the ones in \cite{RO}, \cite{CD5}, \cite{BachmannFasel}, \cite{EHKSY3} and \cite{ElmantoKolderup}, where the details are left to the reader. 
 \footnote{Hint: use the description of the spectrum associated with $\uW_S^\Lambda$. }
\end{rem}

One can formally extend the preceding identification to any base.

\begin{df}\label{df:htp_Witt_singular}
Let $S$ be an arbitrary scheme.
Let $\Lambda \subseteq \QQ$ be a subring with $2 \in \Lambda^\times$.
 We define the $\Lambda$-linear homotopical Witt spectrum $\HW \Lambda_S \in \DA(S,\Lambda)$ as the inverse image $p^*(\HW \Lambda_\ZZ)$ along the projection $p:S \rightarrow \spec(\ZZ)$.
\end{df}

\begin{num}\label{num:homotopical_Witt}
By construction, $\HW \Lambda_S$ defines an absolute ring spectrum (in the sense of Paragraph~\ref{num:abs_pur_def}) as $S$ varies.
Note that Theorem~\ref{thm:Witt motives} implies that the unit map of $\HW \Lambda_S$ induces a canonical isomorphism
\[
  \HW \Lambda_S \simeq \Lambda_{S,-},
\]
and that there is a canonical equivalence
 $$\DA(S,\Lambda)_- \simeq \modl{\HW \Lambda_S},$$
over every base $S$.
In particular, whenever $S$ is regular, we get by Theorem~\ref{thm:Witt motives} the canonical isomorphism
\begin{equation}\label{eq:htp_Witt&regular}
  \HW \Lambda_S \simeq \Lambda_{S,-} \simeq \HuW^\Lambda_S,
\end{equation}
thus recovering Definition~\ref{df:HW}.
\end{num}

\begin{rem}\label{rem:htp_Witt_singular}
For $S$ regular, the identification of Definitions~\ref{df:HW} and \ref{df:htp_Witt_singular} through the isomorphism \eqref{eq:htp_Witt&regular} means in particular that the motivic spectrum $\HW \Lambda_S$ is induced by
 a $\GG$-periodic and strictly $\AA^1$-local sheaf of abelian groups.
 One cannot expect this to hold over singular schemes $S$, due to the failure of the Gersten--Witt conjecture (see Remark~\ref{rem:W_S singular}). In particular, over a singular base $S$ with dimension function $\delta$, it is no longer true that $\HW \Lambda_S$ is concentrated in only one degree for the $\delta$-homotopy t-structure.
\end{rem}
 
\section{Rational stable homotopy and Milnor--Witt motives}
\label{sec:MW}

Based on the initial definition of \cite[Chaps.~3 and 6]{BCDFO} and the comparison result
 of \cite[4.2]{Garkusha}, we introduce the following definition:
\begin{df}\label{df:HMWQ}
Let $S$ be a scheme. One defines the
 \emph{Milnor--Witt rational motivic cohomology spectrum} over $S$, denoted $\HMW \QQ_S \in \SH(S)_\QQ$, by the formula:
$$
\HMW \QQ_S=\HM \QQ_S \oplus \HW \QQ_S
$$
where $\HM \QQ_S$ is the Beilinson motivic cohomology spectrum (cf. \cite[14.1.2]{CD3})
 and $\HW \QQ_S$ is the homotopical Witt spectrum over $S$
 (Definition~\ref{df:htp_Witt_singular}).
This is a motivic ring spectrum, which forms an absolute ring spectrum as $S$ varies.
\end{df}

Recall that $\SH(S)_\QQ$ is canonically identified with $\DA(S,\QQ)$, as symmetric monoidal $\infty$-categories \cite[5.3.35]{CD3}.
An immediate corollary of Theorem~\ref{thm:Witt motives} and \cite[14.2.14]{CD3} is:
\begin{cor}\label{cor:rational_sphere&MW}
For any scheme $S$, the unit map of the ring spectrum $\HMW \QQ_S$ induces an isomorphism of ring spectra in $\SH(S)_\QQ$
$$
\un_S \otimes \QQ \simeq \HMW \QQ_S
$$
compatible with base change.
In particular, for any regular scheme $S$, one obtains canonical isomorphisms
 compatible with pullbacks and products:
$$
\HA^{n,i}(S,\QQ):=\Hom_{\SH(S)_\QQ}(\un_S,\un_S(i)[n])
 \simeq K_{2i-n}^{(i)}(S) \oplus H^{n-i}_\zar(S_\QQ,\uW) \otimes_\ZZ \QQ
$$
for all $(n,i)\in\ZZ^2$.
\end{cor}
Note that this formula was conjectured by Morel (see discussion after Theorem~5.2.2 in \cite{Morelpi0}).
We can deduce formally the following extension of a result due to Garkusha 
 (see \cite{Garkusha}).
\begin{cor}\label{cor:rational_SH&DMW}
For any scheme $S$, there is a canonical equivalence of symmetric monoidal $\infty$-categories:
$$
\SH(S)_\QQ \simeq \modl{\HMW \QQ_S}.
$$
Moreover, the collection of these equivalences for various schemes $S$
 induces an equivalence of motivic $\infty$-categories:
$$
\SH_\QQ \simeq \modl{\HMW \QQ},
$$
which is compatible with the six operations.
\end{cor}

\section{\texorpdfstring{$\SL$}{SL}-Orientations}
\label{sec:SL}

Recall that Panin and Walter have defined several notions that weaken the classical notion of orientation in motivic homotopy theory. We refer the reader to \cite{AnaSL} or \cite[Def.~2.1.3]{DF3}.
In this section, we adapt their definition of $\SL$-orientation to motivic categories, on the model of the usual notion of orientation (see \cite[2.4.38, 2.4.40]{CD3}).

\begin{num}
 We let $\T$ be a motivic $\infty$-category over a category of schemes $\base$.
 Let $S$ be a scheme in $\base$. Recall that the functor $E \mapsto \Th_S(E)$, sending a vector bundle $E$ to its Thom space in $\T(S)$, can be extended to a monoidal functor:
$$
\Th_S:\uK(S) \rightarrow \Pic(\T(S)),
\quad
v \mapsto \Th_S(v).
$$
from the Picard groupoid $\Pic(\T(S)) = \Pic(\T(S),\otimes)$ of virtual vector bundles over $S$ \cite{Del}
 to that of $\otimes$-invertible objects of $\T(S)$.
 Moreover, this assignment is functorial with respect to pullbacks,
 giving a natural transformation $\Th:\uK \rightarrow \Pic(\T)$.
\end{num}

\begin{num}\label{num:Thom&SL_orientation}
An $\SL$-orientation of a virtual vector bundle $v$ over $S$ is an isomorphism
  $$\lambda:\det(v) \simeq \AA^1_S;$$
we refer to the pair $(v,\lambda)$ as an $\SL$-oriented vector bundle.
 The sum of $\SL$-oriented bundles is given by
$$
(v,\lambda)+(v',\lambda')=(v+v',\lambda \otimes v').
$$
The notion of isomorphism of $\SL$-oriented virtual bundles being clear,
 one gets a Picard groupoid $\uK^\SL(S)$ together with a (symmetric) monoidal forgetful
 functor $\psi_S:\uK^\SL(S) \rightarrow \uK(S)$.\footnote{One can also build this category
 from the exact category of $\SL$-torsors.} This is obviously functorial and defines a natural transformation $\psi : \uK^\SL \to \uK$.
\end{num}

Consider the canonical functor:
$$
\twist -_S:\uZZ(S) \rightarrow \Pic(\T(S)),
\quad
r \mapsto \Th_S(\AA^{r}_S)=\un_S\twist r
$$
where $\uZZ(S)$ is the Picard groupoid of locally constant integral functions on $S$.
As $S$ varies this defines a natural transformation $\uZZ \to \Pic(\T)$.

\begin{df}\label{df:SL-orientation}
An \emph{$\SL$-orientation} of the motivic $\infty$-category $\T$ is a natural transformation $\thom$ of presheaves over $\base$
 in monoidal categories, pictured as a $2$-arrow in the following diagram:
$$
\xymatrix@C=20pt@R=-4pt{
\uK^{\SL}\ar^{\Th \circ f}[rrrrdd]\ar_{\rk}[dddd] &&&& \\
& \ar@{=>}^{\thom}[dd] & & & \\
&&&& \Pic(\T) \\
&&&& \\
\uZZ\ar_{\twist -}[rrrruu] &&&&
}
$$
such that for any scheme $S$, $\thom(\AA^1_S)$ is the identity of $\Th(\AA^1_S)$ (where $\AA^1$ is endowed with its tautological orientation).
We refer to the pair $(\T,\thom)$ as an $\SL$-oriented motivic $\infty$-category.
\end{df}

\begin{rem}\label{rem:SL-orientation}
Explicitly, an $\SL$-orientation on $\T$ amounts to the data of, for any scheme $S$ and any $\SL$-oriented virtual bundle $(v,\lambda)$ over $S$,
 a natural isomorphism:
$$
\thom(v,\lambda):\Th_S(v) \rightarrow \un_S\twist{\rk v}.
$$
These isomorphisms are compatible with pullbacks and products.
\end{rem}

\begin{ex}
Clearly a classical orientation of $\T$ induces by restriction
 an $\SL$-orientation.
 In fact, a classical orientation can be formulated just as in Definition~\ref{df:SL-orientation}, except that the upper horizontal arrow is replaced by $\Th: \uK \to \Pic(\T)$.
\end{ex}

\begin{rem}
Recall that the map $\uK \to \Pic(\T)$ lifts to an $\mathcal{E}_\infty$-map from the Thomason--Trobaugh K-theory space $\mathrm{K}$ to the Picard $\infty$-groupoid of $\otimes$-invertible objects in $\T$ (see e.g. \cite[Subsect.~16.2]{BachmannHoyois}).
Let $\mathrm{K}^\SL$ denote the $\SL$-oriented (higher) K-theory presheaf as in \cite[Ex.~16.20]{BachmannHoyois}.
There is again a forgetful map $\mathrm{K}^\SL \to \mathrm{K}$, which we can use to formulate a stronger version of Definition~\ref{df:SL-orientation}.
This amounts to the existence of a homotopy coherent system of compatibilities between the isomorphisms in Remark~\ref{df:SL-orientation}.
We do not discuss this notion further, since we do not know whether the orientation of Theorem~\ref{thm:DA_SL-oriented} is homotopy coherent.
We note however that the orientation on the motivic $\infty$-category of $\MSL$-modules is homotopy coherent, as follows from the motivic Thom spectrum description of $\MSL$ given in \cite[Sect.~16]{BachmannHoyois}.
\end{rem}

\begin{num}
We compare our definition with the original one of Panin and Walter.
 By universality of $\SH$ among motivic $\infty$-categories \cite{Rob}, there is an (essentially unique) adjunction of motivic $\infty$-categories
 \[
 \varphi^*:\SH\leftrightarrows \T:\varphi_*.
 \]
Let $H^{n,i}(S,\T)=\Hom_{\T(S)}(\un_S,\un_S(i)[n])$ be the
 \emph{fundamental cohomology} associated with $\T$.
 For any scheme $S$ in $\base$, $H^{n,i}(-,\T)$ is represented over $S$ by the motivic ring spectrum $\HH^\T_S=\varphi_*(\un_S)\in\SH(S)$.
 The collection $(\HH^\T_S)_{S \in \base}$ form a section of $\SH$ over $\base$.\footnote{Note
 that this section is not cartesian in general. However, it will be in the two cases considered below.}

Then it is easy to check that an $\SL$-orientation $\thom$ of $\T$ is equivalent
 to a family of $\SL$-orientations $\thom^{PW}_S$ of the ring spectra $\HH^\T_S$
 in the sense of Panin--Walter, for all schemes $S$ in $\base$, such that for any morphism of schemes $f:T \rightarrow S$,
 the canonical morphism:
$$
f^*(\HH^\T_S) \rightarrow \HH^\T_T
$$
is compatible with the respective $\SL$-orientations.
 Moreover, if $\base$ admits a final object $\Sigma$,
 then an $\SL$-orientation of $\T$ is uniquely determined by an $\SL$-orientation
 of $\HH^\T_\Sigma$ (by functoriality of $\thom$).
\end{num}

\begin{thm}\label{thm:DA_SL-oriented}
The motivic $\infty$-category $\DA(-,\ZZ)_-$
 admits a unique $\SL$-orientation.
\end{thm}
\begin{proof}
According to the preceding paragraph it is sufficient to give an $\SL$-orientation
 of the ring spectrum $\un_{S,-}$, when $S=\spec{\ZZot}$.
 Then the result follows from \cite[Cor. 5.4]{AnaSL}
 and Corollary \ref{cor:comput_HA1minus}: indeed, according to the latter,
 $H^{0,0}_{\AA^1}(-,\ZZot)_- \simeq \uW^{\ZZot}_S$ is a Zariski sheaf over $\sm_S$.
\end{proof}

As a corollary, we get the analog of the classical fact in topology that rational spectra
 are all orientable.
\begin{cor}\label{cor:SHQ_SL-oriented}
The rational stable homotopy category $\SH_\QQ$ admits a canonical $\SL$-orientation.
\end{cor}
This follows from the decomposition $\SH_\QQ=\DMB \oplus \DA(-,\QQ)_-$.
 Note that the orientation is unique if one asks in addition that its projection on the plus-part
 is induced by the additive orientation on Beilinson motivic cohomology $\HM \QQ$.

\begin{num}
\label{num:grdet}
In the following discussion we consider an arbitrary $\SL$-oriented motivic $\infty$-category $(\T,\thom)$ over $\base$,
 but the reader may focus on the case of $\DA(-,\Lambda)_-$ or $\SH_\QQ$.

Recall that one can define following Deligne \cite{Del} the determinant functor 
$$
\tdet_S:\uK(S) \rightarrow \uPic(S), v \mapsto \big(\det(v),\rk(v)\big)
$$
with values in the Picard groupoid of graded line bundles over $S$. This is again contravariantly functorial in $S$.

Motivated by the following definition, and by the theory of Chow--Witt groups,
 we consider the following twisting functor:
\begin{align}
\begin{split}
\Tw_S:\uPic(S) &\rightarrow \T(S) \\
(L,r) &\mapsto \Th(L \oplus \AA^{r-1}_S)=\Th(L)\twist{r-1}.
\end{split}
\end{align}
We will also write simply $\Tw_S(L)=\Tw_S(L,0)$.

Let $v$ be a virtual vector bundle of rank $r$ over a scheme $S$ in $\base$.
 Then the rank $r-1$ virtual bundle $v-\det v$
 admits an obvious $\SL$-orientation
 $\lambda:\det(v-\det v) \simeq \det v \otimes (\det v)^{-1} \simeq \AA^1_S$.
\end{num}
\begin{df}
Given a virtual vector bundle $v$ over a scheme $S$ in $\base$,
 we define the \emph{twisted Thom class} associated with $\thom$ as the following map $\thom(v) : \Th_S(v) \to \Tw_S(\det v)\langle \rk(v)\rangle$ defined as the composite:
\begin{align*}
\Th_S(v) \simeq \Th_S(\det v) \otimes \Th_S(v-\det v)
 \xrightarrow{1 \otimes \thom(v-\det v,\lambda)} &\Tw_S(\tdet v) \\
 =&\Tw_S(\det v,\rk(v)).
\end{align*}
\end{df}

\begin{rem}\label{rem:weak_orientation}
Note in particular that to give an $\SL$-orientation on $\T$
 in the sense of Definition~\ref{df:SL-orientation}
 is equivalent to give a natural transformation as a 2-arrow in the following diagram:
$$
\xymatrix@C=20pt@R=-4pt{
\uK\ar^{\Th}[rrrrdd]\ar_{\tdet}[dddd] &&&& \\
& \ar@{=>}^{\thom}[dd] & & & \\
&&&& \Pic(\T). \\
&&&& \\
\uPic\ar_{\Tw}[rrrruu] &&&&
}
$$
Note in particular that the twisted Thom class is multiplicative, roughly expressed by the relation
 $\thom(v+w)=\thom(v) \otimes \thom(w)$, using the fact $\Tw_S$ is monoidal (see below).
\end{rem}

\begin{num}
Let $\T$ be an $\SL$-oriented motivic $\infty$-category.
Given line bundles $L$ and $M$ over $S$, there is an induced isomorphism:
$$
\Tw(L) \otimes \Tw(M)=\Th(L \oplus M)\twist{-2}
 \xrightarrow{\thom(L \oplus M)} \Tw(\det(L \oplus M))\simeq \Tw(L \otimes M)
$$
One deduces from this isomorphism that the twisting functor $\Tw_S$ has
 a canonical structure of a monoidal functor:
$$
\Tw_S:\uPic(S) \rightarrow \Pic(\T(S)).
$$
According to a remark of Ananyevskiy (see \cite[Lem.~4.1]{AnaSL}),
 one can go further.
 Indeed, for any line bundle $L$ over $S$,
 one can define a morphism of $S$-schemes
 $\phi_L:L^\times \rightarrow (L^\vee)^\times$ defined locally by sending a non-zero vector $u$
 to the unique function $f$ such that $f(u)=1$.
 The morphism of schemes $\phi_L$ is an isomorphism, so that, by definition of the Thom space, one gets
 a canonical isomorphism in $\T(S)$:
$$
\Th(L) \simeq \Th(L^\vee).
$$
Using the monoidal structure on $\Tw$,
 one gets isomorphisms
$$
\Tw(L) \rightarrow \Tw(L^\vee) \simeq \Tw(L)^{\otimes,-1}
 \Leftrightarrow
 \Tw(L^{\otimes,2}) \simeq \un_S.
$$
In other words, the twisting functor $\Tw_S$ can be factorized and induces a functor from the Picard groupoid
 of $S$ modulo squares.\footnote{This almost shows that the $\SL$-orientation $\thom$ corresponds to an $\SL^c$-orientation
 in the terminology of Panin. See \cite[Rem. 4.4]{AnaSL} for more discussion.}
\end{num}

The following result is obtained by combining the construction of \cite{DJK},
 Theorem \ref{thm:absolute_purity}, Proposition \ref{thm:DA_SL-oriented} and the preceding
 construction.
\begin{thm}\label{thm:SLoriented_absolute_purity}
Let $\T$ be one of the motivic $\infty$-categories $\DA(-,\ZZ)_-$ or $\SH_\QQ$.
Then for any smoothable lci morphism $f:X \rightarrow S$,
 with graded determinant $\tomega_f=\tdet(L_f)$,
 there exists a fundamental class
$$
\eta_f:\Tw_X(\tomega_f) \rightarrow f^!(\un_S)
$$
which is an isomorphism whenever $X$ and $S$ are regular, or $f$ is smooth.
Moreover, the classes $\eta_f$ are compatible with composition, transversal pullbacks, and satisfy an excess intersection formula (see \cite{DJK} for the precise formulations).
\end{thm}

\begin{rem}
As in \cite{DJK}, Theorem~\ref{thm:SLoriented_absolute_purity} gives rise to trace maps and exceptional functoriality (Gysin or wrong-way maps) for the four homology/cohomology theories associated with any spectrum $\E$ in $\DA_-$ or $\SH_\QQ$.
 The only thing that changes is the new twists involved.
 We will illustrate this in the forthcoming section.
\end{rem}

\section{Bivariant \texorpdfstring{$\AA^1$}{A1}-theory and Chow--Witt groups}
\label{sec:bivariant}

\begin{num}
In this section, we will use our results to compute certain $\AA^1$-homology groups.
 We will use two kinds of ring spectra as coefficients,
 which naturally correspond to motivic $\infty$-categories:
\begin{itemize}
\item $R=\HA \ZZ_-$, corresponding to $\T=\DA(-,\ZZ)_-$;
\item $R=\HA \QQ$, corresponding to $\T=\DA(-,\QQ) \simeq \SH_\QQ$.
\end{itemize}
\end{num}
\begin{df}\label{df:bivariant_A1}
Let $R$ and $\T$ be as above.
 Let $X$ be a scheme, $L$ a line bundle over $X$, and $(n,i) \in\ZZ^2$
 a pair of integers.
We define the $L$-th twisted cohomology of $X$ in degree $(n,i)$ and coefficients in $R$ as:
$$
\HA^{n,i}\big(X,R(L)\big)=\Hom_{\T(S)}\big(R_X,R_X\otimes \Tw(L)(i)[n]\big).
$$
Given a morphism
 of $\ZZot$-schemes $f:X \rightarrow S$ that is essentially of finite presentation (Definition~\ref{df:efp}),
 we define the $L$-th twist of the bivariant theory of $X/S$ in degree $(n,i)$ and coefficients in $R$
 as:
$$
\BA_{n,i}\big(X/S,R(L)\big)=\Hom_{\T(S)}\big(R_X\otimes \Tw(L)(i)[n],f^!R_X\big)
$$
where $f^!$ is as defined in Paragraph~\ref{num:f^! efp}.
\end{df}

Taking into account the continuity property of $\T$,
 one gets the following duality theorem as a corollary of Theorem
 \ref{thm:SLoriented_absolute_purity}.\footnote{This formula is obviously reminiscent
 to the duality of Grothendieck--Serre for coherent sheaves.}
\begin{cor}\label{cor:BM=coh}
Let $R$ and $\T$ be as above.
 Then for any smoothable morphism $f:X \rightarrow S$
 between regular schemes,
 of relative dimension $d$ and determinant $\omega_f$,
 any $(n,i) \in \ZZ^2$ and any line bundle $L/X$, one has a canonical isomorphism:
$$
\BA_{n,i}\big(X/S,R(L)\big) \simeq \HA^{2d-n,d-i}\big(X,R(\omega_f-L)\big).
$$
\end{cor}

\begin{num}\label{num:coniveau}
\textit{$\delta$-niveau spectral sequences}.
Now we study the \emph{niveau spectral sequences} in motivic homotopy theory, along the lines of \cite[\S 3.1]{BD1}.
Let us start by fixing our conventions for the next results:
\begin{itemize}
\item All schemes in this section will be assumed noetherian and of finite dimension.
\item $S$ is a regular noetherian scheme and $\delta=-\codim_S$ is the canonical dimension function.
\item $X$ is an essentially of finite type over $S$.
We denote the induced dimension function on $X$ (see \cite[1.1.7]{BD1}) again by $\delta$.
Since it is also assumed noetherian and finite-dimensional, $\delta$ is bounded.
\item \emph{Coefficients}. $R=\HA \QQ$ or $R=\HA \ZZ_-$.
\item \emph{Twists}.
 $L/X$ is a line bundle, $n \in \ZZ$.
\end{itemize}

The $\delta$-niveau spectral sequence of $X/S$ associated with $R_S$ in $\T(S)$
 has the following form:
\begin{equation}\label{eq:nivss}
{}^\delta E^1_{p,q}=\bigoplus_{x \in X_{(p)}} \BA_{p+q,n}\big(\kappa(x)/S,R(L_x)\big) \Rightarrow \BA_{p+q,n}\big(X/S,R(L)\big)
\end{equation}
where $X_{(p)}=\{x \in X \mid \delta(x)=p\}$, for $x \in X$, $\kappa(x)$ is the residue field and $L_x$ is the fiber of $L$ at $x$.
 Given our assumptions, one can apply the previous corollary to 
 the morphism of schemes $\spec{\kappa(x)} \rightarrow S$; note that it is of relative dimension $\delta(x)$,
 and denote by $\omega_x$ its determinant. This gives, for $x \in X_{(p)}$, a canonical isomorphism:
\begin{equation*}
\BA_{p+q,n}\big(\kappa(x)/S,R(L_x)\big) \simeq \HA^{p-q,p-n}\big(\kappa(x),R(\omega_x-L_x)\big)
\end{equation*}
 Taking into account the computation of \eqref{cor:comput_HA1minus}, and the known vanishing of Beilinson motivic cohomology
 of fields, one obtains that the previous spectral sequence is concentrated in one of the regions: ($p>n$, $q>n$)
 or ($p \in [\delta_-(X),\delta_+(X)]$, $q=n$), where $\delta_-(X)$ (resp. $\delta_+(X)$) is the infimum (resp. supremum) of values of $\delta$ over $X$.

Moreover, using Morel's computations of the graded stable homotopy groups of the zero sphere over fields
 as Milnor--Witt K-theory (cf. \cite{MorelA1}),
 one can compute the terms of the complex located on the line $q=n$ as (graded by the index $p$):
\begin{equation*}
\hdots \rightarrow \bigoplus_{x \in X_{(p)}} M_{p-n}\big(x,\omega_x \otimes L_x^\vee\big)
 \xrightarrow{d_{p-1}} \bigoplus_{y \in X_{(p-1)}} M_{p-n-1}\big(y,\omega_y \otimes L_y^\vee\big) \hdots 
\end{equation*}
where we have put, without twists by a line bundle:
\begin{itemize}
\item if $R=\HA \QQ$, $M_r(x):=\KMW_r\big(\kappa(x)\big)_\QQ$;
\item if $R=\HA \ZZ_-$, $M_r(x):=W\big(\kappa(x)\big)[\ot]$;
\end{itemize}
and with twists:
$$
M_r\big(x,\omega_x \otimes L_x^\vee\big):=M_r(x) \otimes_{\ZZ[\kappa(x)^\times]} R[(\omega_x \otimes L_x^\vee)^\times].
$$
\end{num}
\begin{df}
\label{df:CHtilde}
We define the $(L,n)$-th twisted Gersten complex of $X$ with coefficients in $M=\KMW_{*,\QQ}$
 (resp. $M=W[\ot]$),
 with respect to the dimension function $\delta$, as the complex just defined above
 for $R=\HA \QQ$ (resp. $R=\HA \ZZ_-$).
 We denote it by $C^\delta_*\big(X,M(L,n)\big)$,
 and by $A^\delta_p\big(X,M(L,n)\big)$ its homology in degree $p$.
We define the \emph{$\QQ$-linear $n$-th Chow--Witt group of $X$ twisted by $L$ and graded by $\delta$} as:
$$
\wCH_{\delta=n}(X,L)_\QQ:=A^\delta_p\big(X,\KMW_*(L,n)_\QQ\big).
$$
\end{df}

\begin{num}\label{num:RostSchmid}
 Let us consider the case $R=\HA \QQ$.
By construction, the above Gersten complex in degrees around $n$ is of the form:
$$
 \bigoplus_{\eta \in X_{(n+1)}} \KMW_1\big(\eta,\omega_\eta \otimes L_\eta^\vee\big)_\QQ
 \xrightarrow{d_n} \bigoplus_{x \in X_{(n)}} GW\big(x,\omega_x \otimes L_x^\vee\big)_\QQ
 \rightarrow \bigoplus_{s \in X_{(n-1)}} W\big(s,\omega_s \otimes L_s^\vee\big)_\QQ
$$
so that $\wCH_{\delta=n}(X,L)_R$ is a sub-quotient of cycles on $X$ with coefficients
 in quadratic forms and $\delta$-dimension $n$.
 Using the theory of \cite{Feld2},
 one can check that the differentials have the expected form,
 of the so-called \emph{Rost-Schmid complex},\footnote{See \cite[\S 3.3]{Feld2}
 stated in the case where $X=S$ is a smooth scheme over a field.
 The general case works similarly.}
 so that this definition extends the previously known one of \cite{Fasel_wCH} when $X=S$
 is regular. In the case where $S$ is a smooth $k$-scheme, the reader can find the relevant
 result in \cite[Chap. 6, Lemma 4.2.3]{BCDFO}.

Note also that in degree $p>n$, the Gersten complexes with coefficients in
 $(\KMW_*)_\QQ$ and in $W_\QQ$ agree. In the end, using the vanishing
 of the $\delta$-niveau spectral sequence as explained in the paragraph preceding the definition, one gets:
\end{num}
\begin{thm}\label{thm:biv_A1&ChowWitt}
Let $S$ be a regular scheme, $\delta=-\codim_S$ and $X$ an essentially of finite type $S$-scheme.
 Then one gets canonical isomorphisms:
$$
\BA_{r,n}\big(X/S,\QQ(L)\big) \simeq
\begin{cases}
\wCH_{\delta=n}(X,L)_\QQ & r=2n, \\
A^\delta_{n-1}\big(X,\KMW(L,n)_\QQ\big) & r=2n-1, \\
A^\delta_i\big(X,\uW(L)_\QQ\big) & r=n+i, i>n.
\end{cases}
$$
\end{thm}

\begin{rem}\label{rem:Gersten_functoriality}
The $\delta$-niveau spectral filtration
 is well-known to be functorial in the $S$-scheme $X$,
 covariantly with respect to proper maps,
 contravariantly with respect to \'etale maps.
 With a bit more work, as in \cite{Jin},
 one can show that the above $\delta$-niveau spectral sequence, and therefore
 the Gersten complex $C_*^\delta(X,\KMW(L,n))$ satisfies the same functoriality.
 One then obtains that the above isomorphisms are natural with respect to
 that functoriality on the right hand-side and the usual functoriality of
 bivariant theories (as in \cite{DJK} in the motivic case).

With a little more efforts, one can get the contravariant functoriality
 in $X/S$ of the Gersten complex with respect to smoothable lci maps
 as in Feld's theory of Milnor--Witt modules (see \cite{Feld1}).
 The above isomorphisms are again compatible with this functoriality,
 where on the left hand-side one uses the Gysin maps of \cite{DJK}.
\end{rem}

It is interesting to rewrite the case $X=S$.
\begin{cor}
Let $X$ be a regular scheme. Then one has canonical isomorphisms:
\begin{align*}
\Hom_{\SH(X)}(\un_X,\Tw_X(L)(n)[2n]) \otimes_\ZZ\QQ \simeq \wCH^n(X,L)_\QQ \\
 \simeq \CH^n(X)_\QQ \oplus H^n_\zar(X_\QQ,\uW(L))_\QQ
\end{align*}
where $\wCH^n(X,L)$ denotes Chow--Witt groups in codimension $n$ with rational coefficients,
 and $\uW(L)$ denotes the $L$-twisted Witt sheaf.
\end{cor}

Note also that in the case $R=\HA \ZZ_-$ the spectral sequence degenerate at $E_2$
 and gives the following:
\begin{thm}\label{thm:biv_A1minus&Wittcoh}
Let $S$ be a regular scheme, $\delta=-\codim_S$ and $X$ an essentially of finite type $S$-scheme.
 Then one gets a canonical isomorphism:
$$
\BA_{r,n}\big(X/S,\ZZ_-(L)\big) \simeq
A^\delta_{r-n}\big(X_\QQ,W(L,n)[\ot]\big).
$$
\end{thm}

\begin{rem}
The Gersten complex of $X_\QQ$ with coefficients in $W(L,n)[\ot]$ can be identified
 with Gille's version of the Gersten--Witt complex associated with the dualizing sheaf
 $f^{-1}(\omega_S) \otimes L$ as defined in \cite{Gille}.
\end{rem}

\begin{num}
We give here an application to the study of weight structures on motivic categories. Let $\T$ be a motivic $\infty$-category and let $S$ be a scheme. Recall that the \emph{Chow weight structure} on $\T(S)$, if it exists, is the unique weight structure on the underlying motivic triangulated category whose heart agrees with the idempotent completion of the localizing subcategory generated by objects of the form $f_*\un_X(r)[2r]$, where $r\in\ZZ$ and $f:X\to S$ is a proper morphism with $X$ regular (see \cite[Theorem 2.1]{Bondarkorweight}).

We consider the question of the existence of the Chow weight structure on $\T=\SH_\QQ$.
Using de Jong's alterations, if $S$ is separated of finite type over an excellent scheme of dimension $\leqslant 2$, and $-1$ is a sum of squares in any residue field of $S$, then $\SH(S)_{\QQ}\simeq\DM(S)_{\QQ}$ has a Chow weight structure by \cite{Hebert} and \cite{Bondarkorweight}. In what follows we show that the converse is almost true, provided that a general version of resolution of singularities is available.

\begin{lm}
\label{lm:BMWnonvan}
Let $X_0\to S_1$ be a proper, dominant, generically finite morphism between integral schemes of characteristic $0$, with $X_0$ regular. Assume that the function field $\kappa$ of $X_0$ is formally real. Then the group $\BA_{0,0}\big((X_0\times_{S_1}X_0)/X_0,\QQ_-\big)$ is non-vanishing.
\end{lm}
\proof
Denote $T=X_0\times_{S_1}X_0$. By Theorem~\ref{thm:biv_A1minus&Wittcoh}, the group $\BA_{0,0}\big(T/X_0,\QQ_-\big)$ agrees with the kernel of the first differential in the Gersten-Witt complex
$$
\BA_{0,0}\big(T/X_0,\QQ_-\big)
=
\operatorname{ker}\big(\bigoplus_{x\in T^{(0)}}W(\kappa(x))_\QQ\to\bigoplus_{x\in T^{(1)}}W(\kappa(x),\omega_x)_\QQ\big).
$$
The residue fields of the generic points of $T$ are all isomorphic to $\kappa$, and the group $W(\kappa)_\QQ$ is non-vanishing by \cite[VIII 6.15]{LamQuadratic}.
This implies that the group $\BA_{0,0}\big(T/X_0,\QQ_-\big)$ is non-vanishing: indeed, let $\eta$ be a generic point of $T$, then the structure map 
$$\QQ\simeq W(\QQ)_\QQ\to W(\kappa(\eta))_\QQ$$ 
factors through the kernel $\BA_{0,0}\big(T/X_0,\QQ_-\big)$, and for an ordering $\sigma$ of $\kappa(\eta)$ the composition with the signature map on $\kappa(\eta)$
$$\QQ\simeq W(\QQ)_\QQ\to \BA_{0,0}\big(T/X_0,\QQ_-\big)\to\QQ$$
agrees with the signature on $\QQ$ induced by $\sigma$, which shows that $\BA_{0,0}\big(T/X_0,\QQ_-\big)$ is non-zero.
\endproof

\begin{cor}
\label{cor:chowwtnotexist}
Let $S$ be a scheme and let $s$ be a point of $S$. Assume that there exists a proper morphism $f:X\to S$ with $X$ regular, such that
\begin{itemize}
\item the induced map between the characteristic $0$ fibers $X_0\to S_0$ factors through $S_1$, the reduced closure of $s$ in $S_0$;
\item the morphism $X_0\to S_1$ satisfies the condition in Lemma~\ref{lm:BMWnonvan}.
\footnote{Note that these assumptions imply that the residue field of $s$ is formally real.}
\end{itemize}
Then there cannot exist a Chow weight structure on $\SH(S)_{\QQ}$.
\end{cor}
\proof
If the Chow weight structure exists, then the objects $f_*\un_X(r)[2r]$ are negative in the sense of \cite[Definition 1.5 VII]{Bondarkorweight}. In particular this means that the groups $\BA_{2r-n,r}\big((X\times_{S}X)/X,\QQ_-\big)\simeq\BA_{2r-n,r}\big((X_0\times_{S_1}X_0)/X_0,\QQ_-\big)$ vanish for all $n>0$. Since the spectrum $\HA \ZZ_-$ is $\mathbb{G}_m$-periodic by Theorem~\ref{thm:representability_Witt}, we have in particular $\BA_{0,0}\big((X\times_{S}X)/X,\QQ_-\big)=0$, which contradicts Lemma~\ref{lm:BMWnonvan}.
\endproof

\begin{rem}\leavevmode

\begin{enumerate}
\item
The criterion of Corollary~\ref{cor:chowwtnotexist} is to be understood as follows: if $s$ is a formally real point of $S$, let $S_0$ be the characteristic $0$ fiber of $S$ and let $S_1$ be the closure of $s$ in $S_0$. Then Hironaka's resolution of singularities says that there exists a proper birational morphism $X_0\to S_1$ with $X_0$ regular. Corollary~\ref{cor:chowwtnotexist} states that if $X_0$ has a regular model over the whole $S$, then there is no Chow weight structure on $\SH(S)_{\QQ}$.

\item
The condition in Corollary~\ref{cor:chowwtnotexist} holds when $S$ has equal characteristic, or when $S$ is the spectrum of a Dedekind domain (see \cite[1.2.2]{Scholl}). The condition is almost a de Jong alteration, except that we require the morphism to be generically finite of \emph{odd} degree, since a finite field extension of odd degree of a formally real field remains formally real.
\end{enumerate}
\end{rem}

\end{num}

\appendix
\section{Continuity in motivic \texorpdfstring{$\infty$}{∞}-categories}
\label{apd:continuity}

  \begin{setup}\label{setup:continuity}
    Suppose given a cofiltered system $(S_\alpha)_\alpha$ of schemes with affine transition maps.
    We write $S$ for its limit, $u_\alpha : S \to S_\alpha$ for the projection map, and $u_{\alpha,\beta} : S_\beta \to S_\alpha$ for the transition map induced by a morphism $\alpha \to \beta$ in the indexing category. The maps are labeled as follows:
    \[
      \begin{tikzcd}
        S \ar{r}{u_\beta}\ar[swap]{rd}{u_\alpha}
        & S_\beta \ar{d}{u_{\alpha,\beta}}
        \\
        & S_\alpha.
      \end{tikzcd}
    \]
  \end{setup}

  \begin{df}\label{df:continuity}
    Let $\T$ be a motivic $\infty$-category.
    For every cofiltered system $(S_\alpha)_\alpha$ as in Setup~\ref{setup:continuity}, consider the canonical functor
    \[
      \ilim_\alpha \T(S_\alpha) \to \T(S),
    \]
    where the colimit is taken in the $\infty$-category of presentable $\infty$-categories and left-adjoint functors \cite[Def.~5.5.3.1]{LurieHTT}, and the transition maps are the inverse image functors $u_{\alpha,\beta}^*$.
    We say that $\T$ is \emph{continuous} if this functor is an equivalence for all such systems $(S_\alpha)_\alpha$.
  \end{df}

  \begin{num}\label{num:continuity lower star}
    By \cite[Cor.~5.5.3.4]{LurieHTT}, the condition in Def.~\ref{df:continuity} is equivalent to the assertion that the canonical functor
    \[
      \T(S) \to \plim_\alpha \T(S_\alpha)
    \]
    is an equivalence, where the limit is taken in the $\infty$-category of presentable $\infty$-categories and \emph{right}-adjoint functors, and the transition maps are the direct image functors $(u_{\alpha,\beta})_*$.
    By \cite[Thm.~5.5.3.18]{LurieHTT}, this limit can be taken equivalently in the $\infty$-category of (large) $\infty$-categories.

    In particular, the space of maps between any two objects $\E,\F \in \T(S)$ can be computed by the limit
    \[
      \Maps_{\SH(S)}(\E, \F)
      \simeq \plim_\alpha \Maps_{\SH(S_\alpha)}((u_\alpha)_*(\E), (u_\alpha)_*(\F)).
    \]
  \end{num}

  \begin{ex}\label{ex:continuous}
    The motivic $\infty$-category $\SH$ is continuous \cite[Prop.~C.12(4)]{HoyoisLefschetz}.
    We immediately deduce the same for $\SH[1/2]$, $\SH_+$, $\SH_-$, $\SH_\QQ$, $\DA(-,\Lambda)$ (for any commutative ring $\Lambda$), etc.
  \end{ex}

  \begin{df}
    We say that $\T$ is \emph{compactly generated} if the stable $\infty$-category $\T(S)$ is compactly generated for every scheme $S$ \cite[Def.~5.5.7.1]{LurieHTT}, and for any morphism $f : S \to T$, the inverse image functor $f^* : \T(T) \to \T(S)$ preserves compact objects.
    We write $\T(S)^\omega$ for the full subcategory of compact objects in $\T(S)$ for any scheme $S$.
    This is a thick stable subcategory, and we have $\Ind(\T(S)^\omega) \simeq \T(S)$.
  \end{df}

  \begin{rem}\label{rem:f_* colimits}
    Suppose that $\T$ is compactly generated.
    Then by adjunction, the direct image functor $f_* : \T(S) \to \T(T)$ commutes with filtered colimits for any morphism $f : S \to T$.
    Since it is an exact functor between stable $\infty$-categories, it therefore commutes with arbitrary small colimits.
  \end{rem}

  \begin{df}\label{df:constructibly generated}
    Recall from \cite[Def.~4.2.1]{CD3} the full subcategory $\T_{\mathrm{c}}$ of \emph{constructible} objects.
    We say that a motivic $\infty$-category $\T$ is \emph{constructibly generated} if it is compactly generated and the constructible objects are compact.
    Then we have $\T_c(S) = \T(S)^\omega$ for every scheme $S$, i.e., the compact and constructible objects coincide \cite[Prop.~1.4.11]{CD3}.
  \end{df}

  \begin{rem}\label{rem:f_! constructibles}
    Suppose that $\T$ is constructibly generated.
    In this case, the exceptional inverse image functor $f^!$ commutes with colimits for any morphism $f$ of finite type.
    This follows by adjunction from the fact that $f_!$ preserves compact objects \cite[Cor.~4.2.12]{CD3} (where the separatedness and noetherian hypotheses are not necessary, see \cite[Cor.~2.62]{KhanSix}).
  \end{rem}

  \begin{ex}
    The motivic $\infty$-category $\SH$ is constructibly generated by \cite[Prop.~C.12(1-3)]{HoyoisLefschetz}.
    We deduce the same for $\SH[1/2]$, $\SH_+$, $\SH_-$, $\SH_\QQ$, $\DA(-,\Lambda)$ (for any commutative ring $\Lambda$), etc.
  \end{ex}

  \begin{num}\label{num:continuity compact}
    Let $\T$ be a compactly generated motivic $\infty$-category.
    It follows from Prop.~5.5.7.8 and Cor.~4.4.5.21 in \cite{LurieHTT} that $\T$ is continuous if and only if, for every cofiltered system $(S_\alpha)_\alpha$ as in Setup~\ref{setup:continuity}, the canonical functor
    \[
      \ilim_\alpha \T(S_\alpha)^\omega \to \T(S)^\omega
    \]
    is an equivalence, where the colimit is taken in the $\infty$-category of (small) $\infty$-categories and the transition arrows are the inverse image functors $u_{\alpha,\beta}^*$.
    Concretely, we have:
    \begin{enumerate}
      \item\label{item:continuity compact maps}
      For any index $\alpha_0$ and objects $\E_{\alpha_0},\F_{\alpha_0} \in \T(S_{\alpha_0})$ with $\E_{\alpha_0}$ \emph{compact}, we have canonical isomorphisms of spaces
      \begin{equation*}
        \ilim_{\alpha\geqslant\alpha_0} \Maps_{\T(S_\alpha)}(\E_\alpha, \F_\alpha)
        \to \Maps_{\T(S)}(\E, \F),
      \end{equation*}
      where $\E_\alpha \in \T(S_\alpha)$ and $\E \in \T(S)$ are the inverse images (and similarly for $\F_\alpha$ and $\F$).
      (Using the compactness of $\E_{\alpha_0}$, we can reduce to the case where $\F_{\alpha_0}$ is also compact.
      Then the claim holds by the fully faithfulness of the functor above, since formation of mapping spaces in an $\infty$-category commutes with filtered colimits.)

      \item
      For every compact object $\E \in \T(S)^\omega$, there exists an index $\alpha$, a compact object $\F_\alpha \in \T(S_\alpha)^\omega$, and an isomorphism $\E \simeq u_\alpha^*(\F_\alpha)$.
    \end{enumerate}
  \end{num}

  We finish this appendix by formulating continuity in terms of the cohomology and bivariant theories represented by an object of a motivic $\infty$-category.

  \begin{notation}\label{notation:cohomology+bivariant}
    Let $\T$ be a motivic $\infty$-category and $\E$ an absolute object of $\T$.
    Recall the following notation from \cite[Def.~2.2.1]{DJK}: for any scheme $X$ and virtual vector bundle $v$ on $X$, the $v$-twisted cohomology space of $X$ is
    \[
      \E(X, v) = \Maps_{\T(X)}\big( \un_X, \E_X \otimes \Th_X(v) \big),
    \]
    where $\Th_X$ is the Thom space functor in $\T$ (see Paragraph~\ref{num:Thom&SL_orientation}).
    For any finite type morphism $f : X \to S$, we define the $v$-twisted bivariant theory space as
    \[
      \E(X/S, v)
      = \Maps_{\T(X)}\big( \Th_X(v), f^!\E_S \big)
      \simeq \Maps_{\T(S)}\big( f_!\Th_X(v), \E_S \big).
    \]
  \end{notation}

  \begin{num}\label{num:continuity cohomology}
    Let $\T$ be a constructibly generated, continuous motivic $\infty$-category.
    Let $(S_\alpha)_\alpha$ be a cofiltered system of schemes as in Setup~\ref{setup:continuity}, with limit $S$.
    By continuity for algebraic K-theory \cite[Thm.~7.2]{TT}, any virtual vector bundle $v$ on $S$ is the inverse image of a virtual vector bundle $v_{\alpha_0}$ on $S_{\alpha_0}$, for some sufficiently large index $\alpha_0$.
    For any absolute object $\E$ of $\T$, Paragraph~\ref{num:continuity compact}(\ref{item:continuity compact maps}) yields that the canonical maps
    \[
      \ilim_{\alpha\ge\alpha_0} \E(S_{\alpha}, v_{\alpha})
      \to \E(S, v)
    \]
    are invertible (since the unit $\un_{S_{\alpha_0}} \in \T(S_{\alpha_0})$ is constructible and hence compact).
  \end{num}

  \begin{num}\label{num:continuity bivariant}
    Let $\T$ be a constructibly generated, continuous motivic $\infty$-category.
    Let $(S_\alpha)_\alpha$ be a cofiltered system of schemes as in Setup~\ref{setup:continuity}, with limit $S$.
    Let $X_{\alpha_0}$ be a finite type $S_{\alpha_0}$-scheme for some index $\alpha_0$, let $X_\alpha = X_{\alpha_0} \times_{S_{\alpha_0}} S_\alpha$ for every $\alpha\ge\alpha_0$, and let $X \simeq X_{\alpha_0} \times_{S_{\alpha_0}} S$ be the limit of $(X_\alpha)_{\alpha\ge\alpha_0}$.
    Denote by $f_\alpha : X_\alpha \to S_\alpha$ and $f : X \to S$ the structural morphisms.
    Let $v$ be a vector bundle over $X$, which as in Paragraph~\ref{num:continuity cohomology} we may assume (by modifying $\alpha_0$ if necessary) is the inverse image of some virtual vector bundle $v_{\alpha_0}$ on $X_{\alpha_0}$.
    Then the canonical maps
    \[
      \ilim_{\alpha\ge\alpha_0} \E(X_\alpha/S_\alpha, v_\alpha)
      \to \E(X/S, v)
    \]
    are invertible.
    Indeed, the constructible generation assumption on $\T$ implies (see Remark~\ref{rem:f_! constructibles}) that $(f_\alpha)_!(\Th_{X_\alpha}(v_\alpha))$ is compact for every $\alpha$.
    Moreover, by the base change formula, $f_!\Th_X(v)$ is the inverse image of $(f_\alpha)_!(\Th_{X_\alpha}(v_\alpha))$ for every $\alpha$.
    Hence the claim follows from Paragraph~\ref{num:continuity compact}(\ref{item:continuity compact maps}).
  \end{num}

\section{Essentially of finite presentation morphisms}
\label{apd:efp}

  \subsection{Essentially of finite presentation and smooth morphisms}

    The material in this subsection is well-known, but does not seem to be recorded in the literature yet.
    We recall our global conventions that all schemes are implicitly assumed to be quasi-compact and quasi-separated.

    \begin{df}\label{df:pro-etale}
      A morphism of schemes $f : X \to Y$ is \emph{pro-finitely presented} if it can be written as the limit of a cofiltered system of finitely presented morphisms $f_\alpha : X_\alpha \to Y$, with affine transition maps.
      We define \emph{pro-smooth} and \emph{pro-\'etale} morphisms, and \emph{pro-open immersions} analogously.
      Note that the transition maps are automatically \'etale in the latter two cases.
    \end{df}

    \begin{df}\label{df:efp}
      Let $f : X \to Y$ be a morphism of schemes.
      We say that $f$ is \emph{essentially of finite presentation} (resp. \emph{essentially of finite type}) if it factors through a pro-\'etale morphism $p : X \to X'$ and a morphism of finite presentation (resp. of finite type) $g : X' \to Y$.
      We define \emph{essentially smooth} and \emph{essentially \'etale} morphisms analogously.
    \end{df}

    \begin{ex}\label{ex:residue efp}
      If $X$ is a scheme and $x \in X$ is a (not necessarily closed) point, then the inclusion of the residue field $i_x : \spec{\kappa(x)} \to X$ is essentially of finite type via the canonical factorization
      \[
        i_x : \{x\}
        \xrightarrow{j_x} \overline{\{x\}}
        \xrightarrow{k_x} X.
      \]
      If $X$ is noetherian, then $i_x$ is essentially of finite presentation.
    \end{ex}

    \begin{prop}\label{prop:efp compose}
      The class of essentially of finite presentation (resp. essentially smooth, essentially \'etale) morphisms is stable under base change and composition.
    \end{prop}

    \begin{prop}\label{prop:pro-fp=efp}
      A morphism $f : X \to Y$ is essentially of finite presentation (resp. essentially smooth, essentially \'etale) if and only if it can be written as the limit of a cofiltered system of finitely presented morphisms $f_\alpha : X_\alpha \to Y$ (resp. smooth morphisms, \'etale morphisms), with affine and \emph{\'etale} transition maps.
    \end{prop}

    In particular, Proposition~\ref{prop:pro-fp=efp} says that ``essentially \'etale'' is equivalent to ``pro-\'etale''.
    Both Propositions~\ref{prop:efp compose} and \ref{prop:pro-fp=efp} are simple consequences of the following lemma:

    \begin{lm}\label{lm:pro-et composition}
      Let $\mathcal{P}$ be the class of finitely presented morphisms (resp. smooth morphisms, \'etale morphisms, open immersions).
      Let $\mathcal{Q}$ be the class of morphisms $f : X \to Y$ that can be written as the limit of a cofiltered system of morphisms $f_\alpha : X_\alpha \to Y$ belonging to $\mathcal{P}$, with affine and \'etale transition maps.
      Then the class $\mathcal{Q}$ is stable under base change and composition.
    \end{lm}
    \begin{proof}
      Stability under base change is clear.
      Let $f : X \to Y$ and $g : Y \to Z$ be morphisms in $\mathcal{Q}$.
      Write them as limits of cofiltered systems $(f_\alpha : X_\alpha \to Y)_\alpha$ and $(g_\beta : Y_\beta \to Z)_\beta$, respectively.
      Then for every index $\alpha$, there exists by \cite[Thm.~8.8.2(ii)]{EGA4} an index $\beta$, a finitely presented morphism $f_{\alpha,\beta} : X_{\alpha,\beta} \to Y_\beta$ and a $Y$-isomorphism $X_{\alpha,\beta} \times_{Y_\beta} Y \simeq X_\alpha$.
      Then $g\circ f : X \to Z$ is the limit of the finitely presented morphisms $X_{\beta,\alpha} \to Y_\beta \to Z$.
      In the respective cases, $f_{\alpha,\beta}$ is smooth, \'etale, or an open immersion by \cite[Prop.~17.7.8 or Prop.~8.6.3]{EGA4}.
    \end{proof}

    \begin{proof}[Proof of Proposition~\ref{prop:pro-fp=efp}]
      The case of \'etale morphisms follows immediately from the fact that pro-\'etale morphisms are stable under composition (Lemma~\ref{lm:pro-et composition}).
      Let $\mathcal{P}$ be the class of smooth or finitely presented morphisms and let $\mathcal{Q}$ be as in Lemma~\ref{lm:pro-et composition}.
      If $f : X \to Y$ is essentially $\mathcal{P}$, then it is by definition the composite of two morphisms in $\mathcal{Q}$, and hence itself belongs to $\mathcal{Q}$ by Lemma~\ref{lm:pro-et composition}.
      Conversely suppose that $f : X \to Y$ belongs to $\mathcal{Q}$.
      Write $f$ as the limit of a cofiltered system of morphisms $f_\alpha : X_\alpha \to Y$ in $\mathcal{P}$, with affine and \'etale transition maps.
      Then for every $\alpha$, the projection $u_\alpha : X \to X_\alpha$ is pro-\'etale, so the factorization
      \[
        f : X \xrightarrow{u_\alpha} X_\alpha \xrightarrow{f_\alpha} Y
      \]
      exhibits $f$ as essentially $\mathcal{P}$.
    \end{proof}

    \begin{proof}[Proof of Proposition~\ref{prop:efp compose}]
      Combine Proposition~\ref{prop:pro-fp=efp} and Lemma~\ref{lm:pro-et composition}.
    \end{proof}

    \begin{lm}\label{lm:pro-et+fp}
      Let $f : X \to Y$ be a pro-\'etale morphism.
      If $f$ is of finite presentation, then it is \'etale.
    \end{lm}
    \begin{proof}
      It will suffice by definition to show that $f$ is formally \'etale.
      If $f$ is the limit of a cofiltered system of \'etale morphisms $f_\alpha : X_\alpha \to Y$, then the relative cotangent complex $L_{X/Y}$ is the colimit over $\alpha$ of $u_\alpha^*L_{X_\alpha/Y}$, where $u_\alpha : X \to X_\alpha$ is the projection (see e.g. \cite[Tags~08S9, 08R2]{Stacks}).
      Hence $L_{X/Y} \simeq 0$, i.e., $f$ is formally \'etale.
    \end{proof}

    \begin{lm}\label{lm:pro-et 2-of-3}
      Let $f : X \to Y$ be a morphism of $S$-schemes.
      If $X$ and $Y$ are pro-\'etale over $S$, then $f$ is pro-\'etale.
    \end{lm}
    \begin{proof}
      Write $X$ and $Y$ as limits of cofiltered systems $(X_\alpha)_\alpha$ and $(Y_\beta)_\beta$ of \'etale $S$-schemes, respectively.
      For every index $\beta$, there exists by \cite[Prop.~8.14.2]{EGA4} an index $\alpha$ and an $S$-morphism $f_\beta : X_\alpha \to Y_\beta$ through which the composite $X \to Y \to Y_\beta$ factors.
      Now $f_\beta$ is \'etale (as a morphism between \'etale $S$-schemes), and $f : X \to Y$ is the limit of the \'etale maps $X_\alpha \times_{Y_\beta} Y \to Y$ (base changes of $f_\beta$).
    \end{proof}

    \begin{lm}\label{lm:efp 2-of-3}
      Let $f : X \to Y$ be a morphism of $S$-schemes.
      If $X$ is pro-\'etale over $S$ and $Y$ is of finite presentation over $S$, then $f$ is essentially of finite presentation.
    \end{lm}
    \begin{proof}
      Write $X$ as the limit of a cofiltered system of \'etale $S$-schemes $X_\alpha$ (with affine transition maps).
      By \cite[Prop.~8.14.2]{EGA4}, $f : X \to Y$ factors through an $S$-morphism $f_\alpha : X_\alpha \to Y$ for some index $\alpha$.
      Then $f$ is the limit of the finitely presented morphisms $f_\beta : X_\beta \to X_\alpha \to Y$ for all $\beta\geqslant\alpha$.
      Hence it is essentially of finite presentation by Proposition~\ref{prop:pro-fp=efp}.
    \end{proof}

  \subsection{\texorpdfstring{$f^*$}{f\textasciicircum*} for pro-\'etale morphisms}

    \begin{lm}\label{lm:pro-open immersion}
      Let $\T$ be a continuous motivic $\infty$-category.
      Then for every pro-open immersion $j : U \to X$, the functor $j_* : \T(U) \to \T(X)$ is fully faithful.
      Equivalently, the co-unit $j^*j_* \to \Id$ is invertible.
    \end{lm}
    \begin{proof}
      Write $j$ as the limit of a cofiltered system of open immersions $(j_\alpha : U_\alpha \to X)_\alpha$.
      Let $\E$ and $\F$ be objects of $\T(U)$.
      By continuity (Paragraph~\ref{num:continuity lower star}), the canonical map
      \[
        \Maps_{\T(U)}(\E, \F)
        \to \Maps_{\T(X)}(j_*(\E), j_*(\F))
      \]
      is identified with the map
      \[
        \plim_\alpha \Maps_{\T(U_\alpha)}((u_\alpha)_*(\E), (u_\alpha)_*(\F))
        \to \Maps_{\T(X)}(j_*(\E), j_*(\F))
      \]
      induced by the functors $(j_\alpha)_*$.
      Since each $(j_\alpha)_*$ is fully faithful by smooth base change \cite[1.1.20]{CD3}, this map is an isomorphism.
    \end{proof}

    If $\T$ is compactly generated, then Lemma~\ref{lm:pro-open immersion} follows alternatively from the following base change formula:

    \begin{prop}\label{prop:pro-etale base change}
      Let $\T$ be a continuous motivic $\infty$-category.
      Assume that for every morphism of schemes $f$, the direct image $f_*$ commutes with colimits (by Remark~\ref{rem:f_* colimits}, this is automatic if $\T$ is compactly generated).
      Suppose given a cartesian square of schemes
      \[
        \begin{tikzcd}
          X' \ar{r}{g}\ar{d}{p}
            & Y' \ar{d}{q}
          \\
          X \ar{r}{f}
            & Y
        \end{tikzcd}
      \]
      with $f$ and $g$ pro-\'etale.
      Then there is a canonical isomorphism of functors
      \[
        f^*q_* \simeq p_*g^*.
      \]
    \end{prop}
    \begin{proof}
      The natural transformation is defined as the composite
      \[
        f^*q_*
        \xrightarrow{\mathrm{unit}} p_*p^*f^*q_*
        \simeq p_*g^*q^*q_*
        \xrightarrow{\mathrm{counit}} p_*g^*.
      \]
      Write $f$ as the limit of a cofiltered system of \'etale maps $f_\alpha : X_\alpha \to X$, and write $u_\alpha : X \to X_\alpha$ for the projections.
      Denote by $g_\alpha : X'_\alpha\to Y$ and $v_\alpha : X' \to X'_\alpha$ the respective base changes of $f_\alpha$ and $u_\alpha$ along $q : Y' \to Y$.
      By continuity (Paragraph~\ref{num:continuity lower star}), it will suffice to show that the induced natural transformation
      \[
        (u_\alpha)_*f^*q_* \to (u_\alpha)_*p_*g^*
      \]
      is invertible for all $\alpha$.
      Since $f = f_\alpha\circ u_\alpha$ and $g = g_\alpha\circ v_\alpha$, this is identified with
      \[
        (u_\alpha)_*u_\alpha^*(p_\alpha)_*g_\alpha^*
        \simeq (u_\alpha)_*u_\alpha^*f_\alpha^*q_*
        \to (u_\alpha)_*p_*v_\alpha^*g_\alpha^*
      \]
      where we have used the base change isomorphism $f_\alpha^*q_* \simeq (p_\alpha)_*g_\alpha^*$ (since $f_\alpha$ is \'etale).
      It remains therefore to demonstrate the base change formula $u_\alpha^*(p_\alpha)_* \simeq p_*v_\alpha^*$.
      For every morphism $\alpha \to \beta$ in the indexing category, we have the base change isomorphism $u_{\alpha,\beta}^*(p_\beta)_* \simeq (p_\alpha)_*v_{\alpha,\beta}^*$ (since the transition maps $u_{\alpha,\beta} : X_\alpha \to X_\beta$ are \'etale), via which the square
      \[
        \begin{tikzcd}
          \T(X'_\beta) \ar{r}{v_{\alpha,\beta}^*}\ar{d}{(p_\beta)_*}
            & \T(X'_\alpha) \ar{d}{(p_\alpha)_*}
          \\
          \T(X_\beta) \ar{r}{u_{\alpha,\beta}^*}
            & \T(X_\alpha)
        \end{tikzcd}
      \]
      commutes.
      Since $(p_\alpha)_*$ is a left adjoint for each $\alpha$ (as $\T$ is compactly generated, see Remark~\ref{rem:f_* colimits}), this can be regarded as a diagram in the $\infty$-category of presentable $\infty$-categories and left-adjoint functors.
      Thus by continuity, taking colimits in both rows results in a commutative square
      \[
        \begin{tikzcd}
          \T(X'_\alpha) \ar{r}{v_{\alpha}^*}\ar{d}{(p_\alpha)_*}
            & \T(X') \ar{d}{p_*}
          \\
          \T(X_\alpha) \ar{r}{u_{\alpha}^*}
            & \T(X)
        \end{tikzcd}
      \]
      providing the desired base change isomorphism $u_\alpha^*(p_\alpha)_* \simeq p_*v_\alpha^*$.
    \end{proof}

    \begin{prop}\label{prop:Ex^*! pro-etale}
      Let $\T$ be a continuous motivic $\infty$-category.
      Assume that for every morphism of schemes $f$ (resp. of finite presentation), the direct image functor $f_*$ (resp. exceptional inverse image functor $f^!$) commutes with colimits.
      (By Remark~\ref{rem:f_! constructibles}, this is automatic if $\T$ is constructibly generated.)
      Then for every commutative square of schemes
      \[
        \begin{tikzcd}
          X' \ar{r}{g}\ar{d}{p}
            & Y' \ar{d}{q}
          \\
          X \ar{r}{f}
            & Y
        \end{tikzcd}
      \]
      with $f$ and $g$ of finite presentation, there is a canonical morphism of functors
      \[
        Ex^{*!} : p^*f^! \to g^!q^*
      \]
      which is invertible if $q$ is pro-\'etale.
    \end{prop}
    \begin{proof}
      The natural transformation is defined as the composite
      \[
        p^*f^!
        \xrightarrow{\mathrm{unit}} p^*f^!q_*q^*
        \simeq p^*p_*g^!q^*
        \xrightarrow{\mathrm{counit}}
        g^!q^*,
      \]
      where the middle isomorphism is Proposition~\ref{prop:pro-etale base change}.

      Note that the canonical map $a : X' \to X \times_Y Y'$ is pro-\'etale by Lemma~\ref{lm:pro-et 2-of-3}, since the source and target are both pro-\'etale over $X$.
      Since the source and target are of finite presentation over $Y'$, it follows that $a$ is also of finite presentation and hence \'etale by Lemma~\ref{lm:pro-et+fp}.
      In particular, we have a canonical isomorphism $a^* \simeq a^!$.
      Thus we reduce to the case where the square is cartesian.

      Write $q$ as a cofiltered limit of \'etale maps $q_\alpha : Y'_\alpha \to Y$, and denote by $v_\alpha : Y' \to Y'_\alpha$ the projections.
      Consider the cartesian square
      \[
        \begin{tikzcd}
          X' \ar{r}{g}\ar{d}{u_\alpha}
            & Y' \ar{d}{v_\alpha}
          \\
          X'_\alpha \ar{r}{g_\alpha}\ar{d}{p_\alpha}
            & Y'_\alpha \ar{d}{q_\alpha}
          \\
          X \ar{r}{f}
            & Y.
        \end{tikzcd}
      \]
      Under the exchange isomorphism $Ex^{*!} : p_\alpha^*f^! \simeq g_\alpha^!q_\alpha^*$ (since $q_\alpha$ is \'etale), the natural transformation in question is identified with
      \[
        (u_\alpha)_*u_\alpha^*g_\alpha^!q_\alpha^*
        \to (u_\alpha)_*g^!v_\alpha^*q_\alpha^*.
      \]
      Thus it will suffice to show that the exchange transformation $Ex^{*!} : u_\alpha^*g_\alpha^! \to g^!v_\alpha^*$ is invertible.
      Just as in the previous proof, this follows from the exchange isomorphisms $Ex^{*!} : u_{\alpha,\beta}^*g_\beta^! \simeq g_\alpha^!v_{\alpha,\beta}^*$ for all $\alpha \to \beta$ in the indexing category (since the transition maps $u_{\alpha,\beta}$ are \'etale), using continuity and the assumption that the functors $g_\alpha^!$ are left adjoints for all $\alpha$.
    \end{proof}

  \subsection{\texorpdfstring{$f^!$}{f\textasciicircum!} for essentially of finite presentation morphisms}

    Classically, the pair of exceptional adjoint functors $(f_!,f^!)$ is defined under the assumption that $f$ is separated of finite type.
    If $\T$ is a motivic $\infty$-category, then one can use descent to drop the separatedness hypothesis (see \cite[Remark~C.14]{HoyoisLefschetz} or \cite[Appendix~A]{KhanVirtual}).
    That is, one has the exceptional operations $(f_!,f^!)$ for any morphism of finite type, and moreover these satisfy the usual base change and projection formulas.

    In this subsection we show that, under the assumption that $\T$ is continuous, we can further extend $f^!$ to essentially of finite presentation morphisms.
    While the paper was under revision, we learned that similar results, in the context of mixed motives, will appear in work of Benedikt Preis on local Grothendieck--Verdier duality.

    Fix a continuous motivic $\infty$-category $\T$.
    The following definition is inspired by \cite[2.2.12]{BBD}:

    \begin{num}\label{num:f^! efp}
      Let $f : X \to Y$ be an essentially of finite presentation morphism.
      Choose a factorization $f = g \circ p$ with $p$ pro-\'etale and $g$ of finite presentation as in Definition~\ref{df:efp}.
      Then set
      \[
        f^! = p^*g^! : \T(Y) \to \T(X).
      \]
    \end{num}

    \begin{lm}\label{lm:indepf^!}
      The functor $f^!$ is independent of the choice of factorization $f = g \circ p$ (up to homotopy).
    \end{lm}
    \begin{proof}
      Let $f = f' \circ p = f'' \circ q$ be two factorizations, with $p : X \to Y'$ and $q : X \to Y''$ pro-\'etale, $f' : Y' \to Y$ and $f'' : Y'' \to Y$ of finite presentation.
      Then $Y' \times_Y Y'' \to Y$ is finite presentation.
      By Lemma~\ref{lm:efp 2-of-3}, the morphism $(p,q) : X \to Y' \times_Y Y''$ is essentially of finite presentation, since the source is pro-\'etale over $Y'$ and the target is of finite presentation over $Y'$.
      Thus it factors through a pro-\'etale morphism $X \to Y'''$ and a finitely presented morphism $Y''' \to Y' \times_Y Y''$.
      Now $X \to Y''' \to Y$ is a third factorization of $f$ which dominates the chosen ones.
      By symmetry, we may thus assume that there exists a morphism $Y' \to Y''$ fitting in a commutative diagram
      \[
        \begin{tikzcd}
          & Y'\ar{dd}{a}\ar{rd}{f'} &
          \\
          X \ar{ru}{p}\ar{rd}{q} & & Y\\
          & Y'' \ar{ru}{f''} & 
        \end{tikzcd}
      \]
      Regarding the left-hand triangle as a commutative square by inserting the identity $X \to X$, Proposition~\ref{prop:Ex^*! pro-etale} applies to give a canonical isomorphism $p^*a^! \simeq q^*$.
      This then yields the desired isomorphism
      \[
        p^*(f')^!
        \simeq p^*a^!(f'')^!
        \simeq q^*(f'')^!.
      \]
    \end{proof}

    \begin{ex}\label{ex:f^! pro-etale}
      If $f : X \to Y$ is pro-\'etale, then there is an equivalence $f^! \simeq f^*$.
    \end{ex}

    Under mild assumptions on $\T$, we will be able to prove:

    \begin{prop}\label{prop:f^! efp compose}
      Suppose that $\T$ is constructibly generated (Definition~\ref{df:constructibly generated}).
      Let $f : X \to Y$ and $g : Y \to Z$ be essentially of finite presentation morphisms.
      Then there is a canonical isomorphism of functors $(g\circ f)^! \simeq f^!g^!$.
    \end{prop}

    In fact, it suffices to assume that for every morphism of schemes $f$ (resp. of finite presentation), the direct image functor $f_*$ (resp. exceptional inverse image functor $f^!$) commutes with colimits (see Remark~\ref{rem:f_! constructibles}).

    \begin{warn}\label{warn:f_! efp}
      Below we will show that much of the usual behaviour of the operation $f^!$ extends to the essentially of finite presentation case.
      However, we caution the reader that we cannot typically expect $f^!$ to admit a left adjoint $f_!$ when $f$ is not of finite type.
    \end{warn}

    \begin{prop}[Base change]\label{prop:base change efp}
      Assume that for every morphism of schemes $f$, the direct image $f_*$ commutes with colimits.
      Suppose given a cartesian square of schemes
      \[
        \begin{tikzcd}
          X' \ar{r}{g}\ar{d}{p}
            & Y' \ar{d}{q}
          \\
          X \ar{r}{f}
            & Y
        \end{tikzcd}
      \]
      with $f$ and $g$ essentially of finite presentation.
      Then there is a canonical isomorphism of functors
      \[
        p_*g^! \simeq f^!q_*.
      \]
    \end{prop}
    \begin{proof}
      Choose a factorization $f = f' \circ a$ with $f'$ of finite presentation and $a$ pro-\'etale.
      Then the claim follows from the base change formula for finite type morphisms, the identification $a^! \simeq a^*$ (Example~\ref{ex:f^! pro-etale}), and Proposition~\ref{prop:pro-etale base change}.
    \end{proof}

    \begin{prop}[$*!$-exchange]\label{prop:Ex^*! efp}
      Assume that for every morphism of schemes $f$ (resp. of finite presentation), the direct image functor $f_*$ (resp. exceptional inverse image functor $f^!$) commutes with colimits.
      Suppose given a commutative square of schemes
      \[
        \begin{tikzcd}
          X' \ar{r}{g}\ar{d}{p}
            & Y' \ar{d}{q}
          \\
          X \ar{r}{f}
            & Y
        \end{tikzcd}
      \]
      with $f$ and $g$ essentially of finite presentation.
      Then there is a canonical morphism of functors
      \[
        Ex^{*!} : p^*f^! \to g^!q^*
      \]
      which is invertible if $q$ is pro-\'etale.
    \end{prop}
    \begin{proof}
      The natural transformation is defined as the composite
      \[
        p^*f^!
        \xrightarrow{\mathrm{unit}} p^*f^!q_*q^*
        \simeq p^*p_*g^!q^*
        \xrightarrow{\mathrm{counit}}
        g^!q^*,
      \]
      where the middle isomorphism is Proposition~\ref{prop:base change efp}.
      By factoring $f$ through a pro-\'etale morphism and a finite presentation morphism, we immediately reduce to the case where $f$ is of finite presentation, which is Proposition~\ref{prop:Ex^*! pro-etale}.
    \end{proof}

    \begin{proof}[Proof of Proposition~\ref{prop:f^! efp compose}]
      Follows directly from Proposition~\ref{prop:Ex^*! efp}.
    \end{proof}

    We have the following $!$-analogue of \cite[Cor.~4.3.17]{CD3} (with almost the same proof):

    \begin{prop}\label{prop:shriek conservative}
      Let $X$ be a noetherian scheme of finite dimension.
      As $x$ ranges over the points of $X$, the family of functors
      \[
        i_x^! : \T(X) \to \T(\spec{\kappa(x)})
      \]
      is jointly conservative.
    \end{prop}
    \begin{proof}
      Argue by induction on the dimension $d$.
      For $d\leqslant 0$ the claim is clear.
      Assume $d>0$.
      The localization of $X$ at a point $x$ is a pro-open immersion $j_x : \spec{\OO_{X,x}} \to X$.
      Since $j_x^! \simeq j_x^*$ (Example~\ref{ex:f^! pro-etale}), we may use \cite[Prop.~4.3.9]{CD3} (and Proposition~\ref{prop:f^! efp compose}) to reduce to the case where $X$ is local.
      Then the claim follows easily using the localization triangle
      \[
        (i_{x_0})_*i_{x_0}^! \to \Id \to j_*j^!,
      \]
      where $i_{x_0}$ is the inclusion of the closed point (which is essentially of finite presentation, see Example~\ref{ex:residue efp}), and $j$ is the inclusion of the complement (which is of dimension $< d$).
    \end{proof}

\section{Absolute purity in equal characteristic}
\label{apd:popescu}

Let $\T$ be a motivic $\infty$-category fibred over a category of schemes $\base$, and $\E$ an absolute spectrum in $\T$.
As recalled in~\ref{num:abs_pur_def}, for $f:X\to Y$ a smoothable morphism of regular schemes with virtual tangent bundle $T_f$, there is a map
\begin{align}
\label{eq:pur_fdl_two}
\eta_f:\E_X\otimes\Th(T_f)\xrightarrow{}f^!\E_Y
\end{align}
induced by the fundamental class of $f$ \cite[Definition~4.3.7]{DJK}. 
We say that $\E$ is \emph{$f$-pure} if the map~\eqref{eq:pur_fdl_two} is an isomorphism, and we say that $\E$ satisfies \emph{absolute purity} if it is $f$-pure for every smoothable morphism $f$ between regular schemes in $\base$.
In fact, it suffices to check this property for closed immersions \cite[Remark~4.3.12(i)]{DJK}.

This appendix is dedicated to the following result:

\begin{thm}
\label{thm:abs_pur_eqchar}
Let $k$ be a perfect field and let $\T$ be a continuous, constructibly generated motivic $\infty$-category (Defs.~\ref{df:continuity} and \ref{df:constructibly generated}) fibered over $k$-schemes.
Then any absolute spectrum in $\T$ satisfies absolute purity.
\end{thm}

This statement was first asserted in \cite[Example~1.3.4(2)]{DegliseOrientation} (see also \cite[Proposition~6.2]{CD5}).
It is proven, as suggested in \emph{loc. cit.}, by using Popescu's theorem to reduce to relative purity (the case of a closed immersion between \emph{smooth} schemes).
In this appendix we record the details of two variants of this proof.

Recall that a homomorphism of commutative rings $u:A\to B$ is \emph{regular} if it is flat and for every residue field $A \to \kappa$, the ring $B\otimes_A\kappa$ is geometrically regular. Then Popescu's theorem asserts (see \cite[Theorem~1.1]{Spivakovsky}):

\begin{thm}[Popescu]
\label{thm:popescu}
Let $A$ and $B$ be noetherian commutative rings.
A ring homomorphism $u:A \to B$ is regular if and only if it is ind-smooth, i.e., it exhibits $B$ as a filtered colimit of smooth $A$-algebras.
\end{thm}

Before proceeding to the proofs, we begin with a few generalities.

\begin{df}
\label{df:trans_cl_pair}
A \emph{regular closed pair} $(X,Z)$ is a regular closed immersion $i:Z\to X$ between schemes. A \emph{transverse morphism} of regular closed pairs $f:(Y,W)\to (X,Z)$ is a Tor-independent cartesian square
\begin{equation*}
  \begin{tikzcd}
    W \ar{r}\ar{d}
      & Y \ar{d}
    \\
    Z \ar{r}
      & X.
  \end{tikzcd}
\end{equation*}
\end{df}

\begin{lm}
\label{lem:pur_limit}
Let $\T$ be a continuous, constructibly generated motivic $\infty$-category.
Let $\E$ be an absolute object in $\T$.
Let $(X,Z)$ be a regular closed pair which is the limit of a cofiltered system $(X_\alpha,Z_\alpha)_{\alpha\in\Lambda}$ of regular closed pairs with affine and transverse transition morphisms.
If $\E$ is $i_\alpha$-pure for every $\alpha\in\Lambda$, then $\E$ is $i$-pure.
\end{lm}

\begin{proof}[Proof of Lemma~\ref{lem:pur_limit}]
The assumptions on $\T$ imply (see Paragraph~\ref{num:continuity bivariant}) that the canonical map of spaces
\begin{equation*}
  \ilim_{\alpha\in\Lambda} \E(Z_\alpha/X_\alpha, \vb{-N_{i_\alpha}}) \to \E(Z/X, \vb{-N_i})
\end{equation*}
is invertible (see Notation~\ref{notation:cohomology+bivariant}).
By the transverse base change formula \cite[Theorem~3.2.21]{DJK}, each base change morphism $\E(Z_\alpha/X_\alpha, \vb{-N_{i_\alpha}}) \to \E(Z/X, \vb{-N_i})$ sends the fundamental class $\eta_{i_\alpha}$ to the fundamental class $\eta_i$.
Since each of the classes $\eta_{i_\alpha}$ is invertible by assumption, the claim follows.
\end{proof}

\begin{lm}
\label{lm:pur_ess_sm}
Let $\T$ be a continuous, constructibly generated motivic $\infty$-category.
Let $\E$ be an absolute object in $\T$.
If $i$ is a closed immersion between schemes that are essentially smooth over a scheme $S$, then $\E$ is $i$-pure.
\end{lm}
\proof
The case where both schemes are smooth over $S$ follows from relative purity as in \cite[Prop.~4.3.10(i)]{DJK}.
In general, write $X$ and $Z$ as the limits of projective systems of smooth $S$-schemes $X_\alpha$ and $Z_\alpha$, respectively, with affine \'etale transition morphisms.
Then by \cite[Th\'eor\`eme~8.8.2(ii), Proposition~8.6.3]{EGA4}, the closed immersion $i$ descends to a closed immersion $i_\lambda : Z_\lambda \to X_\lambda$ for some sufficiently large $\lambda$.
In particular, the pair $(X, Z)$ can be written as the limit of a projective system of closed pairs of smooth $S$-schemes, with affine \'etale transition morphisms.
Thus the claim follows from Lemma~\ref{lem:pur_limit} (note that \'etale morphisms of pairs are transverse).
\endproof

\begin{num}
We now give the first proof.
We will require the following lemma:

\begin{lm}
\label{lm:local_gen}
Let $A \to B$ be a local homomorphism of regular local rings.
Let $I\subset A$ be a proper ideal such that $B/IB$ is regular, and set $n = \dim(B) - \dim(B/IB)$.
If there exist $n$ elements $f_1,\ldots,f_n \in A$ generating $I$, then the $f_i$ form a subset of a regular system of parameters for $A$, and the ring $A/I$ is regular.
\end{lm}

\begin{proof}
Let $\mathfrak{m}$ and $\mathfrak{n}$ be the respective maximal ideals of $A$ and $B$.
By assumption, the ideal $IB$ is generated by the images $g_i \in B$ of $f_i$, and $\dim(B/IB) = \dim(B) - n$.
By \cite[0, Proposition~16.3.7]{EGA4}, we deduce that the $g_i$ form a subset of a system of parameters for $B$.
Since $B/IB$ is regular, this implies that they are $B/\mathfrak{n}$-linearly independent as elements of $\mathfrak{n}/\mathfrak{n}^2$ \cite[0, Proposition~17.1.7]{EGA4}.
This implies that the $f_i$ are $A/\mathfrak{m}$-linearly independent in $\mathfrak{m}/\mathfrak{m}^2$.
The claim now follows by another application of \cite[0, Proposition~17.1.7]{EGA4}.
\end{proof}

\begin{prop}
\label{prop:aff_reg_pur}
Let $k$ be a perfect field and let $i:Z\to X$ be a closed immersion between affine regular $k$-schemes. Then the closed pair $(X,Z)$ is a projective limit of closed pairs of smooth $k$-schemes with transverse affine transition morphisms.
\end{prop}
\proof
Let $X=\operatorname{Spec} A$ and let $I$ be the ideal defining $Z$. By Theorem~\ref{thm:popescu}, we may write $A$ as a colimit of smooth $k$-algebras $A_\alpha$ for $\alpha\in \Lambda$.
Truncating $\Lambda$ if necessary, we may assume that for every $\alpha\in\Lambda$, there is an ideal $I_\alpha$ of $A_\alpha$ such that $I=I_\alpha A$.
Let $X_\alpha=\operatorname{Spec} A_\alpha$ and $Z_\alpha=\operatorname{Spec} A_\alpha/I_\alpha$.

Let $n$ be the codimension of $Z$ in $X$, considered as a Zariski-locally constant function on $Z$.
Truncating $\Lambda$ again, we may also assume that for all $\alpha\in\Lambda$, the ideal $I_\alpha$ is locally generated by at most $n$ elements. Indeed, the locus $T_\alpha$ where this is not true is closed in $Z_\alpha$. As $Z$ has codimension $n$ in $X$, the projective limit of $T_\alpha$ is empty, so we have $T_\alpha=\emptyset$ for all sufficiently large $\alpha$.

Now let $X'_\alpha$ be the open subscheme of $X_\alpha$ given by the complement of the singular locus of $Z_\alpha$, and let $Z'_\alpha=X'_\alpha\cap Z_\alpha$. It follows from Lemma~\ref{lm:local_gen} that
\begin{itemize}
\item $X'_\alpha$ contains the image of $X$ in $X_\alpha$, and consequently the map $X\to X_\alpha$ factors through $X'_\alpha$.
\item The induced map $X\to X'_\alpha$ is a transverse morphism of regular closed pairs $(X,Z)\to (X'_\alpha,Z'_\alpha)$, since both pairs are defined by regular closed immersions of the same codimension. The transition morphisms of the projective system $(X'_\alpha,Z'_\alpha)_{\alpha\in\Lambda}$ are also transverse for the same reason.
\end{itemize}
Since $k$ is perfect, $(X'_\alpha,Z'_\alpha)$ is a closed pair of smooth $k$-schemes. Furthermore, since $(X,Z)$ is the projective limit of $(X_\alpha,Z_\alpha)_{\alpha\in\Lambda}$, it is also the projective limit of $(X'_\alpha,Z'_\alpha)_{\alpha\in\Lambda}$, and the result follows. 
\endproof

\begin{proof}[First proof of Theorem~\ref{thm:abs_pur_eqchar}]
  Let $\T$ be a motivic $\infty$-category as in the statement.
  By Zariski descent \cite[Proposition~3.3.4]{CD3}, it suffices to show that $\E$ is $i$-pure for any closed immersion $i:Z\to X$ between affine regular $k$-schemes. This follows by combining Lemma~\ref{lem:pur_limit}, Lemma~\ref{lm:pur_ess_sm} and Proposition~\ref{prop:aff_reg_pur}.
\end{proof}
\end{num}

\begin{num}
The second proof is similar but uses an argument of Gabber to first reduce to the case of punctual purity \cite[XVI 3.2]{Gabber}.

\begin{lm}\label{lm:local_gen_max}
Let $A \to B$ be a local homomorphism of regular local rings.
If there exist $n$ elements $f_1,\ldots,f_n \in A$ whose images $g_1,\ldots,g_n \in B$ form a regular sequence of parameters for $B$, then $(f_1,\ldots,f_n)$ is a regular system of parameters for $A$.
In particular, we have $\dim(A) = \dim(B) = n$.
\end{lm}

\begin{proof}
This is a special case of Lemma~\ref{lm:local_gen}, but it is simpler to argue directly as follows.
Since $g_1,\ldots,g_n$ form a minimal set of generators for the maximal ideal $\mathfrak{n} \subset B$, they induce a basis of the $B/\mathfrak{n}$-vector space $\mathfrak{n}/\mathfrak{n}^2$.
It follows that the $f_i$ induce a basis of the $A/\mathfrak{m}$-vector space $\mathfrak{m}/\mathfrak{m}^2$.
By \cite[0, Proposition~17.1.7]{EGA4} this means that the $f_i$ form a regular system of parameters for $A$, as claimed.
\end{proof}

We use the following variant of Proposition~\ref{prop:aff_reg_pur}:
\begin{prop}
\label{prop:local_pur}
Let $A$ be a regular local $k$-algebra with maximal ideal $\mathfrak{m}$. Then the regular closed pair $(\operatorname{Spec} A, \operatorname{Spec} A/\mathfrak{m})$ is a projective limit of closed pairs of essentially smooth affine $k$-schemes with transverse transition morphisms.
\end{prop}
\proof
By Theorem~\ref{thm:popescu}, we may write $A$ as the colimit of an inductive system $(A^0_\gamma)_{\gamma\in\Gamma}$ of smooth $k$-algebras.
If $A_\alpha$ are the localizations of the $A^0_\gamma$'s at prime ideals, then $A$ is also the colimit of the inductive system $(A_\alpha)_{\alpha\in\Lambda}$ of essentially smooth local $k$-algebras.
Let $f_1,\ldots,f_n \in \mathfrak{m}$ be a regular system of parameters for $A$, where $n = \dim(A)$.
For some sufficiently large $\lambda\in\Lambda$, the $f_i$ belong to $A_\lambda$ and form a regular system of parameters of length $n$, by Lemma~\ref{lm:local_gen}.
In particular, $\mathfrak{m} = \mathfrak{m}_\alpha A$ for each $\alpha\geqslant\lambda$, where $\mathfrak{m}_\alpha$ is the maximal ideal of $A_\alpha$, and the closed pair $(\spec{A}, \spec{A}/\mathfrak{m})$ is the limit of the projective system
\begin{equation*}
  (\spec{A_\alpha}, \spec{A_\alpha}/\mathfrak{m}_\alpha)_{\alpha\geqslant\Lambda}.
\end{equation*}
All pairs under consideration are given by regular closed immersions of codimension $n$, so it follows that the transition morphisms are transverse.
\endproof

\begin{proof}[Second proof of Theorem~\ref{thm:abs_pur_eqchar}]
Let $i:Z\to X$ be a closed immersion between regular $k$-schemes. Let $z$ be a point of $Z$, and let $Z_z=\operatorname{Spec}\mathcal{O}_{Z,z}$ and $X_z=\operatorname{Spec}\mathcal{O}_{X,i(z)}$. Denote the canonical closed immersions
\begin{align}
\begin{gathered}
  \xymatrix@=10pt{
    z \ar[r]^-{i_Z} \ar[rd]_-{i_X} & Z_z \ar[d]^-{i'}\\
     & X_z.
  }
\end{gathered}
\end{align}
by $i_Z$ and $i_X$ respectively.
By \cite[Lemma~4.3.17]{CD3} and \cite[Proposition~4.3.9]{CD3}, we only need to show that $i_Z^*\eta_{i'}$ is an isomorphism. By Lemma~\ref{lem:pur_limit}, Lemma~\ref{lm:pur_ess_sm} and Proposition~\ref{prop:local_pur}, we know that both $\eta_{i_Z}$ and $\eta_{i_X}$ are isomorphisms. By \cite[Theorem 3.3.2]{DJK} the fundamental classes are compatible with composition, and we know that $i_Z^!\eta_{i'}$ is an isomorphism.

Now use induction on the dimension of $Z_z$. Let $j_Z:U_z\to Z_z$ be the open immersion given by the complement of $z$ in $Z_z$. Then $U_z$ has dimension lower than $Z_z$, and the induction hypothesis implies that $j_Z^*\eta_{i'}$ is an isomorphism. From the canonical localization triangle
\begin{align}
i_Z^!\eta_{i'}\xrightarrow{} i_Z^*\eta_{i'}\xrightarrow{} i_Z^*j_{Z*}j^*_Z\eta_{i'}\xrightarrow{}i_Z^!\eta_{i'}[1]
\end{align}
and the five lemma, we know that $i_Z^*\eta_{i'}$ is an isomorphism, which finishes the proof.
\end{proof}
\end{num}

\section{Absolute purity in hermitian K-theory}
\label{apd:hermitianK}

In this appendix, we prove absolute purity for hermitian K-theory.
 Recall that Panin and Walter have constructed in \cite{PaninWalter1} for any regular $\ZZot$-scheme $S$
 a ring spectrum $\KQ_S$ which represents Schlichting's higher Grothendieck--Witt groups (\emph{i.e.} hermitian K-theory)
 in $\SH(S)$.
 According to \cite{SchlichTri}, this ring spectrum has a geometric model
 (given by orthogonal grassmanians, see \textit{op. cit.}, Introduction, Th. 1).
 This implies the $\KQ_S$ for various regular $\ZZot$-schemes $S$ are compatible with base change.
 Therefore, one can define $\KQ_S$ for general $\ZZot$-scheme $S$ just by pulling back 
 $\KQ_{\ZZot}$ along the canonical map $S \rightarrow \spec(\ZZot)$.
 Then $\KQ$ is an absolute ring spectrum with respect to $\SH$
 over the category of $\ZZot$-schemes.

\begin{thm}\label{thm:abs_purity_KQ}
The absolute ring spectrum $\KQ$ satisfies absolute purity over $\ZZot$-schemes (in the sense recalled in Paragraph \ref{num:abs_pur_def}).
\end{thm}

The proof follows the same lines as the analogous result for algebraic K-theory \cite[Thm.~13.6.3]{CD3}, and boils down to d\'evissage in Grothendieck--Witt theory.

\begin{proof}[Proof of Theorem~\ref{thm:abs_purity_KQ}]
Let $i:Z\to X$ be a closed immersion of pure codimension $c$ between regular schemes.
The claim is that the morphism $\eta_i : \KQ_Z \otimes \Th_Z(-N_i) \to i^!(\KQ_X)$ defined by the fundamental class is invertible.
Since $i_*$ is conservative and the stable $\infty$-category $\SH(X)$ is generated by objects of the form $\Sigma^\infty_X(T)_+(n)[2n]$, where $T$ varies among smooth $X$-schemes and $n\in\ZZ$, it will suffice to show that the induced map
\begin{equation*}
  [\Sigma^\infty_Z(T_Z)_+(n)[2n],\KQ_Z\otimes\Th_Z(-N_i)]
    \to [\Sigma^\infty_Z(T_Z)_+(n)[2n],i^!(\KQ_X)]
\end{equation*}
is invertible, where $T_Z = T\times_X Z$.
In the notation of \cite{DJK}, this is the map
\begin{equation*}
  \KQ^0(T_Z, -n-\vb{N_i}) \to \KQ^0_{T_Z}(T, -n)
\end{equation*}
induced by pairing with the fundamental class $\eta_{T_Z/Z} \in \KQ_0(T_Z/T, -\vb{N_i})$.
By \cite[Corollary~7.3]{PaninWalter1} and since $\KQ$ is $\SL^c$-oriented, we have canonical isomorphisms
\begin{align*}
  \KQ^0(T_Z, -n-\vb{N_i}) &\simeq \GW_{0}^{[-n-c]}(T_Z, \det(-N_i))\\
  \KQ^0_{T_Z}(T, -n) &\simeq \GW_{0}^{[-n]}(T~\text{on}~T_Z)
\end{align*}
under which the map in question is identified with the Gysin map in Grothendieck--Witt theory induced by direct image of coherent sheaves (see \cite[(9.9)]{Schlicht2}).
The latter is invertible by d\'evissage in Grothendieck--Witt theory \cite[Theorems~9.5, 9.18 and 9.19]{Schlicht2}.
\end{proof}

\section{\texorpdfstring{$2$}{2}-primary unramified Witt ring}
\label{apd:W}

This appendix brings some complements to the theory of Witt and unramified Witt sheaves
 (see Definition~\ref{df:unramified_Witt}) of independent interest. 

\begin{num}\label{num:Balmer}
The theory of Witt groups has been extended to
 higher Witt groups $W^i(S)$ by Balmer (see e.g. \cite{BalmerBasic}).
 Recall that these groups are defined by considering the derived category
 $\T=\Der^b(\mathrm{VB}_S)$ of the exact category of vector bundles over $S$
 equipped with the $(-1)^i$-duality $D_i=\uHom(-,\mathcal O_S)[i]$
 (see \cite[1.4.12]{BalmerBasic}).
 It is then possible to mimic the classical definition of 
 symmetric bilinear forms/symmetric objects, with respect to the duality $D_i$, and metabolic ones (called ``neutral'' in \emph{loc.cit.}).
 One defines $W^i(X)$ as the monoid of isomorphisms classes of
 $D_i$-symmetric objects of $\T$ modulo the submonoid of $D_i$-metabolic ones (see \emph{loc. cit.} 1.4.3).

 The following lemma is the Witt-version of our \emph{key lemma} (Lemma~\ref{lm:fiber_Q_eq}).
\end{num}
\begin{prop}\label{prop:Witt_2inverted}
Let $S$ be a $\ZZot$-scheme and $\Lambda \subseteq \QQ$ a subring with $2 \in \Lambda^\times$.
Let $\nu:S_\QQ \rightarrow S$ be the pro-open immersion of the characteristic $0$ part of $S$
 as in Paragraph \ref{num:zero_part}. 

Then the map $\uW_S^\Lambda \rightarrow \nu_*(\uW_{S_\QQ}^\Lambda)$ from \eqref{eq:funct1_uW} is an isomorphism.
\end{prop}
\begin{proof}
It will be sufficient to prove that for any regular scheme $S$, the following map on Witt groups
 is an isomorphism:
$$
W(S)[2^{-1}] \xrightarrow{\nu^*} W(S_\QQ)[2^{-1}].
$$
We will use Balmer-Witt groups with support, as defined in \cite[Sec. 1]{BalmerWitt}.
 Given a closed pair $(X,Z)$, we write $W^i(X,Z)$ for the $i$-th Witt group of $X$ with support in $Z$.

 Consider a closed subscheme $Z$ of $\spec(\ZZ)$, with open complement $U$.
 Write $S_Z$, $S_U$ for the pullbacks of $S/\spec(\ZZ)$ to $Z$, $U$.
 Using the long exact sequence of Witt groups with support \cite[1.7]{BalmerWitt},
 we get an exact sequence:
$$
W^0(S,S_Z) \rightarrow W^0(S) \rightarrow W^0(S_U) \rightarrow W^1(S,S_Z).
$$
Taking the filtered inductive limit over the category of closed subschemes of $S$,
 and using additivity in $S_Z$, we get an exact sequence:
$$
\oplus_{p\neq 2} W^0(S,S_p) \rightarrow W^0(S) \xrightarrow{\nu^*} W^0(S_\QQ) \rightarrow \oplus_{p\neq 2} W^1(S,S_p)
$$
where $p$ runs over the odd primes.

Therefore, to get the lemma, it is sufficient to prove that for any odd prime $p$,
 the triangular Witt group $W^i(S,S_p)[2^{-1}]$ is zero.

Let $\delta$ be the dimension function on $S$ such that $\delta=-\codim_S$
 (see e.g. \cite[\textsection 1.1]{BD1}).
 From the long exact sequence of Witt groups with support (\emph{op. cit.}, Th. 1.6),
 one can build out of the $\delta$-niveau filtration on $S_p$,
 as in \cite[\textsection 3.1]{BD1}, and deduce a converging spectral
 sequence\footnote{This construction is nothing else than the version with support of the Gersten--Witt spectral sequence,
 whose first version can be found in \cite{BalmerWalter}.} of the form:
$$
E^1_{i,j}=\oplus_{x \in S_{p,(i)}} W^{-i-j}(S_{(x)},\{x\}) \Rightarrow W^{-i-j}(S,S_p)
$$ 
where $S_{p,{(i)}}$ denotes the set of points of $x$ such that $\delta(x)=i$.
 Since $S$ is regular and defined over $\ZZot$, purity (d\'evissage) for the Witt groups \cite[Thm.~6.1]{BalmerWalter} yields:
$$
W^{-i-j}(S_{(x)},\{x\})=
\begin{cases}
W(\kappa(x)) & i+j\equiv 0 \mod 4 \\
0 & \text{otherwise.}
\end{cases}
$$
This group is $2$-primary torsion as $\kappa(x)$ has positive characteristic
 (\cite[VIII 3.2 and 6.15]{LamQuadratic}).
In particular $E^1_{i,j}[\ot]$ vanishes and one can conclude. 
\end{proof}

\begin{rem}
Note that the Witt ring $W(k)$ of a characteristic $2$ field is always a $2$-torsion ring.
 In particular, the only obstruction to the above proof to work for an arbitrary regular scheme $S$
 is that Balmer's theory needs $2$ being invertible on $\mathcal O_S$.
\end{rem}

Let us draw the following corollary:
\begin{cor}\label{cor:Gersten--Witt_conj}
Let $S$ be a regular $\ZZot$-scheme and $\Lambda \subseteq \QQ$
 such that $2 \in \Lambda^\times$.
\begin{enumerate}
\item If $S$ is in addition semi-local, 
 the Gersten--Witt complex of $S$ (see \cite{BalmerWalter}) is a resolution of $W(S)$.
\item Points (1) and (2) of Theorem~\ref{thm:representability_Witt}
 hold for $S$ (\emph{i.e.} without assuming the existence of a base field).
\end{enumerate}
\end{cor}
This immediately follows from the previous proposition,
 as it allows one to reduce to the case of a $\QQ$-scheme and then applies
 \cite[4.3]{BalmerWitt} for the first point and Theorem~\ref{thm:representability_Witt}
 for the second one.

\begin{rem}\label{rem:Jacobson}
The Gersten conjecture for the $2$-primary Witt ring of excellent $\ZZot$-schemes
 was already proven by Jacobson (see \cite[5.3]{Jacobson}).
 Our method gives an alternative (somewhat trivial!) proof
 and also extends the result as we do not need $S$ to be excellent.
\end{rem}
 
\bibliographystyle{amsalpha}
\bibliography{purity}

\end{document}